\documentclass[a4paper]{article} %
\usepackage[margin=2.5cm]{geometry} 
\usepackage{multirow}
\usepackage{multicol}
\usepackage{amsmath}
\usepackage{amsthm}
\usepackage{nicefrac}

\usepackage{amssymb}
\usepackage{graphicx}
\usepackage{epstopdf}
\newcommand{\R}{\mathbb{R}}

\usepackage[utf8]{inputenc}
\usepackage{hyperref}
\hypersetup{
	unicode,
	colorlinks,
	urlcolor=blue, 
	linkcolor=blue, 
	citecolor=blue,
	pdfproducer={LaTeX},
	pdfcreator={pdflatex}
}

\newcommand{\tilt}{\theta}   

\newcommand{\bC}{\mathbf{C}}
\newcommand{\bd}{\mathbf{d}}
\newcommand{\bD}{\mathbf{D}}
\newcommand{\bI}{\mathbf{I}}
\newcommand{\bL}{\mathbf{L}}
\newcommand{\bS}{\mathbf{S}}
\newcommand{\be}{\mathbf{e}}
\newcommand{\by}{\mathbf{y}}
\newcommand{\bG}{\mathbf{G}}
\newcommand{\bR}{\mathbf{R}}
\newcommand{\bK}{\mathbf{K}}
\newcommand{\bm}{\mathbf{m}}
\newcommand{\bp}{\mathbf{p}}
\newcommand{\bq}{\mathbf{q}}

\newcommand{\btheta}{\boldsymbol\theta}
\newcommand{\bgamma}{\boldsymbol\gamma}
\newcommand{\bSigma}{\boldsymbol\Sigma}

\newcommand{\bbR}{\mathbb{R}}

\newcommand{\ex}{_{\text{\scriptsize{true}}}}

\title{Automatic nonstationary anisotropic Tikhonov regularization\\through bilevel optimization}
\author{Silvia Gazzola \thanks{Department of Mathematical Sciences, University of Bath, Bath, UK. (\texttt{sg968@bath.ac.uk})} \and Ali Gholami\thanks{Institute of Geophysics, Polish Academy of Sciences, Warsaw, Poland (\texttt{agholami@igf.edu.pl})}}

\begin{document}
	\maketitle
 
\begin{abstract}
Regularization techniques are necessary to compute meaningful solutions to discrete ill-posed inverse problems. 
The well-known 2-norm Tikhonov regularization method equipped with a discretization of the gradient operator as regularization operator penalizes large gradient components of the solution to overcome instabilities. However, this method is homogeneous, i.e., it does not take into account the orientation of the regularized solution and therefore tends to smooth the desired structures, textures and discontinuities, which often contain important information. If the local orientation field of the solution is known, a possible way to overcome this issue is to implement local anisotropic regularization by penalizing weighted directional derivatives. In this paper, considering problems that are inherently two-dimensional, we propose to automatically and simultaneously recover the regularized solution and the local orientation parameters (used to define the anisotropic regularization term) by solving a bilevel optimization problem. Specifically, the lower level problem is Tikhonov regularization equipped with local anisotropic regularization, while the objective function of the upper level problem encodes some natural assumptions about the local orientation parameters and the Tikhonov regularization parameter.  
Application of the proposed algorithm to a variety of inverse problems in imaging (such as denoising, deblurring, tomography and Dix inversion), with both real and synthetic data, shows its effectiveness and robustness.
\end{abstract}

%
%
%
%
%

\section{Introduction}

Ill-posed inverse problems often arise in various areas of applied science and engineering (see, e.g. \cite{Hansen_1998_RDD,Aster_2004_PEI, Tarantola_2005_IPT}). In this paper we are mostly interested in the solution of large-scale ill-conditioned linear systems of equations of the form
\begin{equation}\label{main_eq}
\bd =\bG\bm\ex + \be, 
\end{equation}
coming from the discretization of continuous inverse problems (see, e.g., \cite{Hansen_1998_RDD,Aster_2004_PEI, Tarantola_2005_IPT}). Here, \linebreak[4]$\bd,\,\be\in\R^{M}$ represent known measured data and unknown white Gaussian noise, respectively; the unknown $\bm\ex\in\R^N$ represents a quantity of interest (sometimes also referred to as model); 
$\bG$ is the discrete forward operator linking the unknown to to the data space. We consider formulation \eqref{main_eq} associated to two-dimensional (2D) problems in space, whereby  
the vector $\bm\ex$ is obtained by stacking the columns of a rectangular $N_z\times N_x$ array, such that its length is $N=N_zN_x$.

Due to the presence of noise and the ill-conditioning of the system in \eqref{main_eq}, suitable regularization is required to produce a reasonable estimate of $\bm\ex$  that fits the data through the forward operator, and is consistent with prior information. Tikhonov regularization in general form computes an approximation $\bm(\mu)$ of $\bm\ex$ as
\begin{equation}\label{tikh_gen}
\bm(\mu):=\arg\min_{\bm\in\bbR^N} \frac12\|\bG\bm-\bd\|_2^2+ \frac{\mu}{2}\|\bK\bm\|_2^2,
\end{equation}
where the penalty parameter $\mu\geq 0$ trades off the fit-to-data and regularization terms,  $\|\cdot\|_2$ denotes the vector $2$-norm, and  $\bK\in\bbR^{P\times N}$ is a regularization matrix encoding some prior information about $\bm\ex$. 
The key to a successful regularizer is the inclusion of correct (prior) information about the solution to be computed. In particular, when the latter is known to contain clear anisotropies (e.g., line-like structures), such information can be used to design an appropriate anisotropic regularizer. Namely, if the local (i.e., entry-wise) orientation parameters $\btheta\in\R^N$ of $\bm\ex$ are known, then they can be encoded in parameter-dependent regularizers within Tikhonov regularization in general form \eqref{tikh_gen}, i.e., one should solve
\begin{equation}\label{tikh_genpar}
\bm(\btheta,\mu):=\arg\min_{\bm\in\bbR^N} \frac12\|\bG\bm-\bd\|_2^2+ \frac{\mu}{2}\|\bK(\btheta)\bm\|_{\bSigma}^2,
\end{equation}
where $\bSigma$ is some weighting (diagonal, positive definite) matrix so that, for a vector $\by$, $\|\by\|_{\bSigma}^2=\by^T\bSigma\by$. The vector $\bgamma:=[\btheta^T, \mu]^T$ may be referred to as hyperparameters (see, e.g., \cite{Fessler}), also considering the relation to the hyperparameters arising within a Bayesian approach to the regularization of \eqref{main_eq} (see, e.g., \cite{Stuart}). However in this paper, following the terminology originally introduced in \cite{OID}, we refer to $\bgamma$ as inversion design parameters, or, briefly, inversion parameters, meaning that they fully characterize the inversion method. One can still apply anisotropic regularization when the orientation parameters of $\bm\ex$ are not known a priori, if a strategy to estimate them is provided. 

There is a rich body of literature on anisotropic regularization. Most of it considers non-smooth, sparsity-enforcing regularizers, such as anisotropic TV. In this paper we only consider smooth regularizers expressed with respect to some (weighted) 2-norm. The upside of this approach is that it is intrinsically simpler and cheaper than dealing with 1-norms, but the results are much improved with respect to those obtained adopting Tikhonov without any anisotropic regularization; see also \cite{L2}. Recently, the authors of \cite{Calatroni19} propose a new space-variant anisotropic regularization term for variational image restoration, based on the statistical assumption that the gradients of the target image distribute locally according to a bivariate generalized Gaussian distribution, coupled with a maximum-likelihood-based estimation to set all the parameters automatically, and within an alternating direction method of multipliers (ADMM) framework. 
Many anisotropic regularization methods are derived for specific geophysical applications. In this setting, the estimated orientation parameters can also be used as a side-product for other purposes (such as velocity analysis, interpolation of spatial gaps in the data, and fault detection); see, e.g., \cite{Fomel_2007_VIT}. Such approaches are mainly targetted to denoising and are typically two-step methods, whereby structure tensors are computed directly from the data and structure-oriented filters are then designed to be applied within a smoother; see, e.g., \cite{Hale_2009_SOS,Fomel_2002_AOP}. 

In this paper we develop a new bilevel optimization method for automatically and simultaneously estimating the local orientation field and the regularized solution of a given discrete inverse problem. The lower level functional is the Tikhonov regularized problem, whose regularization term involves local weighted directional derivatives and whose solution is an argument of the upper level objective function. The latter is the sum of three terms, each designed according to the following assumptions: 
(1) the norm of the directional derivative of the regularized solution is minimal when computed along the true local signal orientation; (2) 
the local orientation varies smoothly; (3) a proper value of the Tikhonov regularization parameter can be set according to the discrepancy principle.
The upper level problem also incorporates box constraints for the orientation parameters and nonnegativity constraints for the Tikhonov regularization parameter. Since both the upper and lower level functionals are smooth, the resulting bilevel optimization problem is solved using L-BFGS-B. 
%
Bilevel optimization methods have recently gained popularity for the solution of inverse problems. They usually constitute the backbones of data-driven approaches to regularization, whereby optimal inversion parameters appearing in a variational regularization method are computed by minimizing an error measure with respect to a set of available training data. Such methods are usually referred to as `supervised bilevel learning'; see \cite{de2016bilevel, kunisch2013bilevel, calatroni2017bilevel, antil2020bilevel}. In contrast to such approaches, the bilevel optimization problem established in this paper does not require training data. To the best of our knowledge, there are no other algorithms in the literature that formulate anisotropic regularization as an unsupervised bilevel optimization problem. 
%
 
The remaining part of the paper is organized as follows. In Section \ref{sect: anisoTikh} we present the nonstationary anisotropic Tikhonov regularization considered in this paper. In Section \ref{sect: bilevel} we introduce the new bilevel optimization approach for simultaneous estimation of the local orientation parameters and the regularized solution. The performance of the new solver is assessed in Section \ref{sect: experiments}, where a variety of synthetic test problems in image deblurring, computed tomography and geophysical inversion are considered; the new method is also compared against isotropic Tikhonov regularization. Section \ref{sect:end} proposes some concluding remarks and highlights possible extensions.

\section{Local anisotropic Tikhonov regularization}\label{sect: anisoTikh}

Before defining the local anisotropic Tikhonov regularization method considered in this paper, it is convenient to introduce its isotropic counterpart (further specifying \eqref{tikh_gen}), i.e., the general Tikhonov regularization method 
\begin{equation}\label{genTikh}
\bm(\mu)=
\arg\min_{\bm\in\bbR^N} \frac12\|\bG\bm-\bd\|_2^2+ \frac{\mu}{2}\|\nabla \bm\|_2^2=
\arg\min_{\bm\in\bbR^N} \frac12\|\bG\bm-\bd\|_2^2+ \frac{\mu}{2}\sum_{i=1}^N\|(\nabla \bm)_i\|_2^2,
\end{equation}
where $\nabla$ is a scaled discrete gradient operator that can be expressed as
\begin{equation}\label{eq:gradient}
\nabla=\left[
\begin{array}{c}
\nabla_{\!x}\\
\nabla_{\!z}
\end{array}\right]=\left[\begin{array}{ccc}
\bD_x &\!\!\!\!\!\otimes\!\!\!\!\!& \bI\\
\bI &\!\!\!\!\!\otimes\!\!\!\!\!&\bD_z
\end{array}
\right],
\end{equation}
where $\nabla_{\!x}$ and $\nabla_{\!z}$ are finite difference operators discretizing the first-order horizontal ($x$) and vertical ($z$) derivatives, $\bI$ is the identity matrix (whose size should be clear from the context) and $\bD_y\in\R^{N_y\times N_y}$, $y=x,z$ are scaled finite difference approximations of the first derivative operator in 1D (specifically, these are bidiagonal matrices , with 1 on the diagonal, -1 on the first superdiagonal,  and a row of zeros at the bottom). The regularization term in \eqref{genTikh} enforces smoothness in the gradient of the solution uniformly in any direction since 
\[
 \left\|\bR(\theta_i)\left[
\begin{array}{c}
(\nabla_{\!x}\bm)_i\\
(\nabla_{\!z}\bm)_i
\end{array}\right]\right\|_2=\|\bR(\theta_i)(\nabla\bm)_i\|_2 = \|(\nabla\bm)_i\|_2
\]
for every rotation matrix $\bR(\theta_i)$ defining directional derivatives with respect to an angle $\theta_i$; see also equation \eqref{rotation}.

In this paper, to overcome this issue, we consider the following anisotropic Tikhonov regularization method, obtained by rotating and scaling the gradient vector in the regularization term in \eqref{genTikh}: \begin{equation} \label{AnisoReg}
 \bm(\btheta,\mu)=\arg\min_{\bm\in\bbR^N}  \frac12\|\bG\bm-\bd\|_2^2+ \frac{\mu}{2}\sum_{i=1}^N\|\bR(\theta_i)(\nabla \bm)_i\|_{\bSigma_i}^2,
\end{equation}
where
\begin{equation}\label{rotation}
\bR(\theta_i)=\left[
\begin{array}{cc}
~~\cos(\theta_i) & \sin(\theta_i)\\
-\sin(\theta_i) & \cos(\theta_i)
\end{array}\right]:=\left[
\begin{array}{c}
\bold{x}'(\theta_i)\\
\bold{z}'(\theta_i)
\end{array}\right],
\quad
\bSigma_i=\left[
\begin{array}{cc}
\sigma^{x'}_i & 0\\
0 & \sigma^{z'}_i
\end{array}\right], \quad
 \quad i=1,...,N.
\end{equation}
In the above equation and elsewhere in this paper, $\theta_i\in(-\pi/2,\pi/2]$, $i=1,\dots,N$, are the local image orientation parameters (or tilt angles), which are measured clockwise about the origin. 
$\bR(\theta_i)$ is a rotation matrix that converts any coordinates in the original $(x,z)$ system into the rotated system $(x',z')$ given in terms of ${\tilt}_i$. Without loss of generality we assume that the entries of the diagonal matrix $\bSigma_i$ are such that $\sigma^{x'}_i\geq \sigma^{z'}_i \geq 0$. This weighting allows for flexible smoothing of the signal along the orientation angle $\theta_i$, as well as control over the anisotropic behavior of the regularization. Indeed, by taking $\sigma^{z'}_i=0$, no smoothing orthogonal to the ${\tilt}_i$ direction is applied and the regularization is maximally anisotropic, so to favor models with elongated features along $\theta_i$. In a Bayesian framework, such regularizer corresponds to assuming that the gradients of the solution are locally distributed according to a bivariate generalized Gaussian distribution; see \cite{Calatroni19}. Although, in general, $\sigma^{x'}_i$ and $\sigma^{z'}_i$ can be adjusted locally for each model point, in this paper we use fixed values; namely, we set 
\begin{equation}\label{weights}
\sigma^{x'}_i=\sigma^{x'}=1\quad\mbox{and}\quad \sigma^{z'}_i=\sigma^{z'}=\epsilon,\quad\mbox{with $\epsilon \ll 1$},    
\end{equation}
forcing line-like linear structures in the final model. 

We conclude this session by noting that problem \eqref{AnisoReg} can be equivalently and compactly rewritten in the form \eqref{tikh_genpar}, namely, as the problem of finding $\bm$ such that
\begin{equation} \label{m_est}
\left(\bG^T\bG+{\mu}\bD_{\epsilon}(\btheta)^T\bD_{\epsilon}(\btheta)\right)\bm=
\bG^T\bd,
\end{equation}
where, using the assumption \eqref{weights} for the weights,  
\begin{equation} \label{D_aniso}
\bD_{\epsilon}(\btheta) = 
\left[\begin{array}{cc}
\bC(\btheta) & \bS(\btheta)\\
-\sqrt{\epsilon} \bS(\btheta) & \sqrt{\epsilon}\bC(\btheta)
\end{array}\right]
\left[\begin{array}{c}
\nabla_x\\
\nabla_z
\end{array}
\right]\,,
\end{equation}
and where $\bC(\btheta)$ and $\bS(\btheta)$ are diagonal matrices with $\cos(\btheta)$ and $\sin(\btheta)$ on their main diagonal, respectively. This alternative expression may be conveniently used to compute the lower level minimizer via either a direct or iterative linear system solver. 

\section{Bilevel optimization\\for model and inversion parameters recovery}\label{sect: bilevel}
To simultaneously reconstruct $\bm$ and recover the inversion parameters $\bgamma=[\btheta^T, \mu]^T$ that define the regularization scheme \eqref{AnisoReg}, we solve a bilevel optimization problem of the form
\begin{equation}\label{eq:bilevel}
\bgamma^\ast=\arg\min_{\bgamma\in\mathcal{C}} U(\bm^\ast(\bgamma),\bgamma)\quad\mbox{subject to}\quad
\bm^\ast(\bgamma)=\arg\min_{\bm\in\mathbb{R}^N}L(\bm,\bgamma)\,,
\end{equation}
In the following we define the quantities appearing in the above equation, motivating their choices. The upper level constraint set $\mathcal{C}=\{\bgamma\in\mathbb{R}^{N+1}\;|\; -\pi/2\leq \gamma_i\leq\pi/2,\: i=1,\dots,N,\:\gamma_{N+1}\geq 0\}$ bounds the orientation parameters and the regularization parameter. The upper level objective function is given by 
\begin{equation}\label{ULobjfcn}
U(\bm^\ast(\bgamma),\bgamma) = \frac 12 \left(\left(\|\bG \bm^\ast(\bgamma)-\bd\|_2^2 - \varepsilon^2\right)^2+\delta^2\right)^{\nicefrac{1}{2}}+ \frac{\alpha}{2} \sum_{i=1}^N\|\bR(\theta_i)({\nabla} \bm)_i\|_{\bSigma_i}^2 + \frac{\beta}{2}\|\nabla \btheta\|_2^2\,.
\end{equation}
The first term in $U(\bm^\ast(\bgamma),\bgamma)$, where $\varepsilon$ is an upper bound for the norm of the noise $\be$ affecting the data $\bd$, enforces a smoothed version of the discrepancy principle (i.e., if the smoothing parameter $\delta=0$, and $\alpha=\beta=0$, then minimizing $U(\bm^\ast(\bgamma),\bgamma)$ would be equivalent to prescribing $\|\bG\bm^\ast(\bgamma)-\bd\|_2^2 = \varepsilon^2$). The second term in $U(\bm^\ast(\bgamma),\bgamma)$ characterizes the orientation angles $\btheta$, as one can show that $\|\bR(\theta_i)({\nabla} \bm)_i\|_{\bSigma}^2$ is minimized when $({\nabla} \bm)_i$ is aligned to $\theta_i$. Indeed, by letting $\bp:=[\cos(\theta_i), \sin(\theta_i)]^T$ and $\bq:=[ -\sin(\theta_i), \cos(\theta_i)]^T$ be the transposes of the rows of $\bR(\theta_i)$, we have
\[
\|\bR(\theta_i)({\nabla}\bm)_i\|_{\bSigma}^2 = (\sigma^{x'} \langle \bp, ({\nabla}\bm)_i\rangle)^2+(\sigma^{z'} \langle \bq, ({\nabla}\bm)_i\rangle)^2,
\]
where $\langle \cdot, \cdot\rangle$ denotes the Euclidean scalar product between vectors. Denoting by $\zeta$ the angle between the vectors $\bp$ and $({\nabla}\bm)_i$ and using well-known properties of trigonometric functions,  the above quantity can also be expressed as
\[
\|\bR(\theta_i)({\nabla}\bm)_i\|_{\bSigma}^2 = ((\sigma^{x'}\cos(\zeta))^2 + (\sigma^{z'} \sin(\zeta))^2)\|({\nabla}\bm)_i\|_2^2\,.
\]
Simple derivations show that the above expression is minimized when $\zeta=0$, i.e., when $({\nabla}\bm)_i$ is aligned to $\bp$. 
The third term in $U(\bm^\ast(\bgamma),\bgamma)$ enforces some smoothness in $\nabla\btheta$, i.e., it penalizes abrupt changes in $\btheta$ by encouraging a constant behavior. Finally, the lower level objective function in \eqref{eq:bilevel} is the objective function of the minimization problem in \eqref{AnisoReg}.

To solve the bilevel optimization problem \eqref{eq:bilevel}, given that the lower level problem is quadratic and has a closed-form solution, we transform it into a single-level optimization problem and apply a gradient-based scheme (L-BFGS) to the resulting problem. We note that
\begin{equation}\label{eq:ULgrad}
\nabla U(\bgamma) = \nabla_{\bgamma}U(\bgamma,\bm^\ast(\bgamma)) + (\nabla_{\bgamma}\bm^\ast(\bgamma))^T\nabla_{\bm}U(\bgamma,\bm^\ast(\bgamma))\,.
\end{equation}
Let $\bL(\gamma_i):=\bSigma^{\nicefrac{1}{2}}\bR(\theta_i)\widehat{\bD}_i$, where
\[
\widehat{\bD}_i=\left[\begin{array}{c}
     \be_i^T(\bD_x\otimes\bI) \\
     \be_i^T(\bI\otimes\bD_z)
\end{array}
\right]\,.
\]
In the above equation $\be_i$ is the $i$th canonical basis vector of $\bR^N$, and $\bI$ and $\bD$ are as in \eqref{eq:gradient}. Then the first term in the above equation is the vector of length $N+1$ whose $i$th component is given by 
\[
(\bm^\ast(\bgamma))^T\left(\left(\frac{\partial \bL(\gamma_i)}{\partial\gamma_i}\right)^T\bL(\gamma_i)+\left(\bL(\gamma_i)\right)^T\frac{\partial \bL(\gamma_i)}{\partial\gamma_i}\right)\bm^\ast(\bgamma)+[\nabla^T\nabla \btheta]_i,\quad \mbox{for $i=1,\dots,N$},
\]
and 0 for $i=N+1$; in the above equation, $\nabla$ is the discrete gradient operator appearing in equation \eqref{ULobjfcn}. The second term in \eqref{eq:ULgrad} is more cumbersome to express, with the $i$th column of the Jacobian $\nabla_{\bgamma}\bm^\ast(\bgamma)$, $i=1,\dots,N$, given by
\begin{equation}\label{grad1}
\frac{\partial \bm^{\ast}(\bgamma)}{\partial\gamma_i}=
-\left(\bG^T\bG+\mu\sum_{i=1}^N\widehat{\bD}_i^T\bR(\theta_i)^T\bSigma_i\bR(\theta_i)\widehat{\bD}_i\right)^{-1}
\mu\left(\left(\frac{\partial \bL(\gamma_i)}{\partial\gamma_i}\right)^T\bL(\gamma_i)+\left(\bL(\gamma_i)\right)^T\frac{\partial \bL(\gamma_i)}{\partial\gamma_i}\right)\,.
\end{equation}
The $(N+1)$th column of the Jacobian is given by 
\begin{equation}\label{grad2}
\frac{\partial \bm^{\ast}(\bgamma)}{\partial\gamma_{N+1}}=
-\left(\bG^T\bG+\mu\sum_{i=1}^N\widehat{\bD}_i^T\bR(\theta_i)^T\bSigma_i\bR(\theta_i)\widehat{\bD}_i\right)^{-1}\left(\sum_{i=1}^N\widehat{\bD}_i^T\bR(\theta_i)^T\bSigma_i\bR(\theta_i)\widehat{\bD}_i\right)\bm^\ast(\bgamma)\,.
\end{equation}
The second factor of the second term in \eqref{eq:ULgrad} is given by
\begin{align*}
\left(\left(\|\bG\bm^\ast(\bgamma)-\bd\|_2^2 - \varepsilon^2\right)^2+\delta^2\right)^{-\nicefrac{1}{2}}\left(\|\bG\bm^\ast(\bgamma)-\bd\|_2^2 - \varepsilon^2\right)\bG^T(\bG\bm^\ast(\bgamma)-\bd)\\
+\sum_{i=1}^N\widehat{\bD}_i^T\bR(\theta_i)^T\bSigma_i\bR(\theta_i)\widehat{\bD}_i\bm^\ast(\gamma)\,.
\end{align*}


We conclude this section by explaining a strategy to avoid instabilities in the estimates of the orientation parameters $\btheta$, which arise when applying the discrete gradient $\nabla$ to the quantity of interest $\bm$ in \eqref{ULobjfcn}.
Specifically, we smooth the gradient components of $\bm$ in the horizontal and vertical directions using appropriate filters ${h}_x$ and ${h}_z$, respectively. Specifically, we define
\begin{equation} \label{dxz}
\widetilde{\nabla}_x\bm={h}_x\ast \bm\quad\mbox{and}\quad \widetilde{\nabla}_z\bm={h}_z\ast \bm,
\end{equation} 
respectively, 
where $\ast$ denotes 2D convolution operator and where
\begin{figure}
\center
\includegraphics[width=9cm]{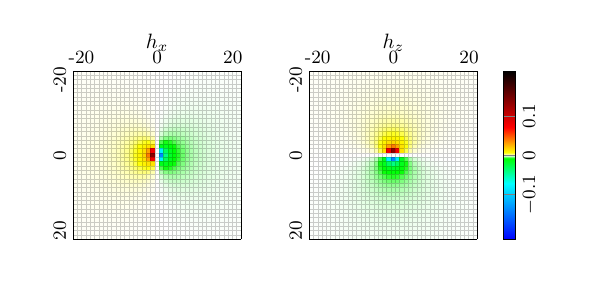}
\caption{The 2D Hilbert transform filters.
}
\label{Riesz_filters}
\end{figure}
\begin{equation}\label{hilbert}
h_x(x,z) = -\frac{1}{2\pi}\frac{x}{(x^2+z^2)^{3/2}},\quad
h_z(x,z) = -\frac{1}{2\pi}\frac{z}{(x^2+z^2)^{3/2}}.
\end{equation}
The above operation is known as Hilbert transform; see also \cite{Claerbout_1976_FGD}. Such smoothed discrete gradient $\widetilde{\nabla}\bm=[(\widetilde{\nabla}_x\bm)^T,(\widetilde{\nabla}_z\bm)^T]^T$ replaces the standard discrete gradient $\nabla\bm$ in \eqref{ULobjfcn} (and, consequently, in the computations of the upper level functional gradient   \eqref{eq:ULgrad}). Different choices of the filters $h_x$ and $h_z$ lead to different smoothed gradients, with the standard gradient (discretized by first-order finite difference) corresponding to
\begin{equation} \label{diffs}
\bold{h}_x=
\begin{bmatrix}
\frac{1}{2} & 0& -\frac{1}{2}
\end{bmatrix},
\quad
\bold{h}_z=
\begin{bmatrix}
\frac{1}{2} &
0&
 -\frac{1}{2}
\end{bmatrix}^T.
\end{equation}
All the operators of the form \eqref{dxz} apply $90^{\circ}$ phase shift and frequency-dependent amplitude scaling that amplifies high frequencies compared to lower frequencies. While the phase shift is the main feature that allows to determine orientation parameters, amplitude scaling causes instabilities when performing this task. Thus, to obtain a stable estimate of the signal variations, we remove the frequency-dependent amplitude scaling of the derivatives and form $\widetilde{\nabla}_x\bm$ and $\widetilde{\nabla}_z\bm$ (with $h_x$ and $h_z$ defined as in \eqref{hilbert}) from the phase information only, which may be regarded as a smoothed derivative. 
Gaussian derivatives are an alternative to the Hilbert transform \eqref{hilbert} in \eqref{dxz}; see \cite{Hale_2006_RGF}. Figure \ref{Riesz_filters} shows the filters in \eqref{hilbert}, were we can clearly see that they are smoothed version of the differentiators in \eqref{diffs}.





\section{Numerical examples}\label{sect: experiments}
In this section we test the performance of the proposed nonstationary anisotropic regularization using several examples.
All the numerical tests are generated using IR Tools \cite{IRtools}. For all experiments, when evaluating the performance of the proposed method, our main term of comparison is a bilevel method, whereby a Tikhonov-regularized problem with a regularization term given by the 2-norm squared of the gradient of the unknown is solved, and the regularization parameter set according to the discrepancy principle. Specifically, the upper level objective function is a function of the Tikhonov regularization parameter, and evaluates the smoothed modulus of the difference between the discrepancy computed for a given regularization parameter and the estimate of the magnitude of the noise; the lower level problem is Tikhonov regularization with a fixed regularization parameter. In all the experiments, $\delta = 0.001$ is the smoothing parameter for the discrepancy principle.

\subsection{Denoising}
\begin{center}
\begin{figure}
\begin{tabular}{cccc}
\hspace{-1.3cm}\small{(a)} & 
\hspace{-1.3cm}\small{(b)} & 
\hspace{-1.3cm}\small{(c)} & 
\hspace{-1.3cm}\small{(d)} \vspace{-0.1cm}\\
\hspace{-1.3cm}\includegraphics[width = 5cm]{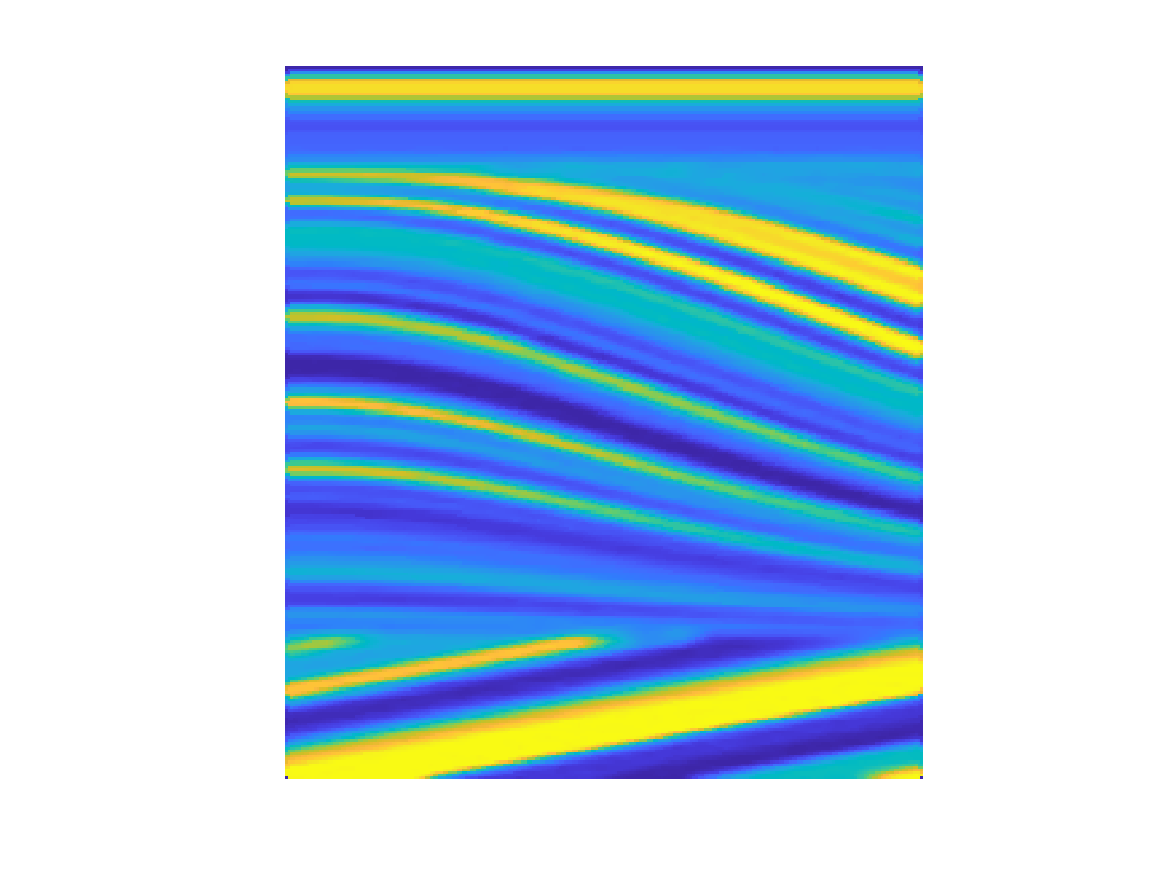} &
\hspace{-1.3cm}\includegraphics[width = 5cm]{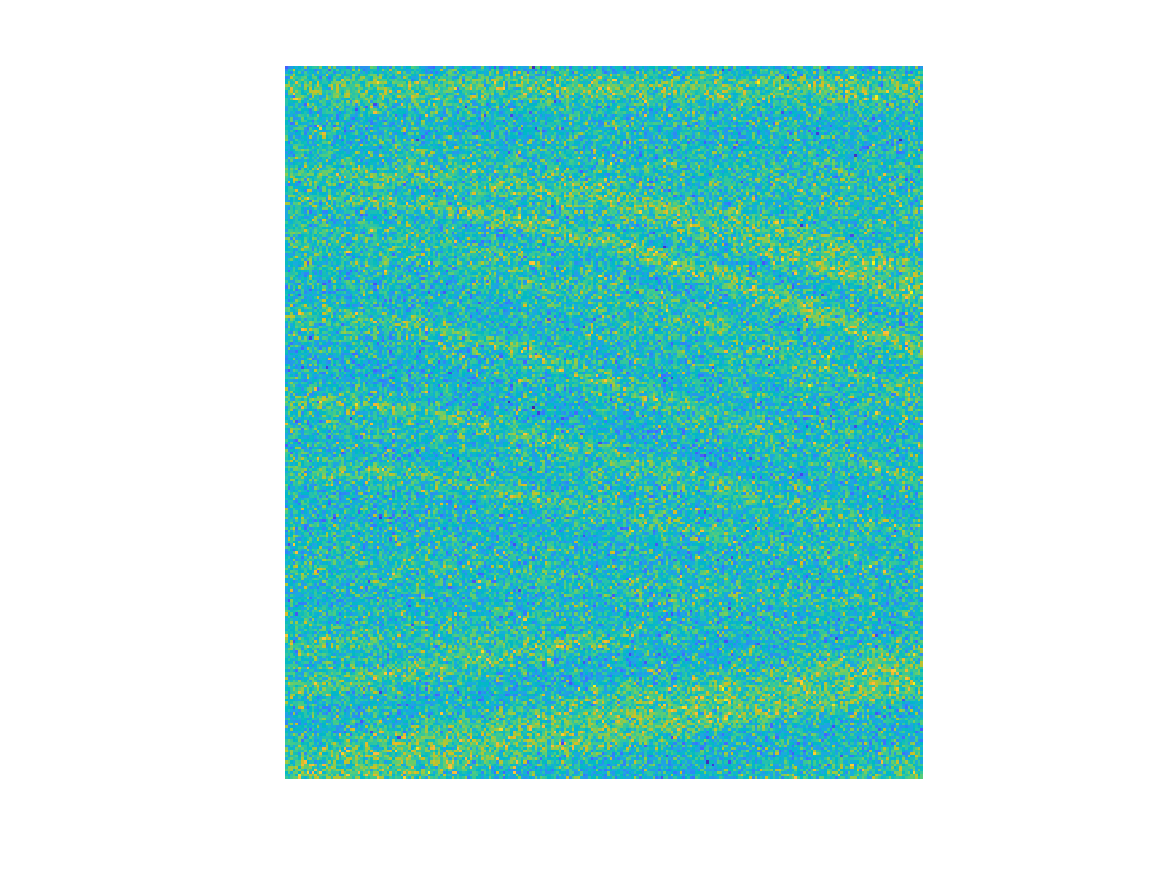} & 
\hspace{-1.3cm}\includegraphics[width = 5cm]{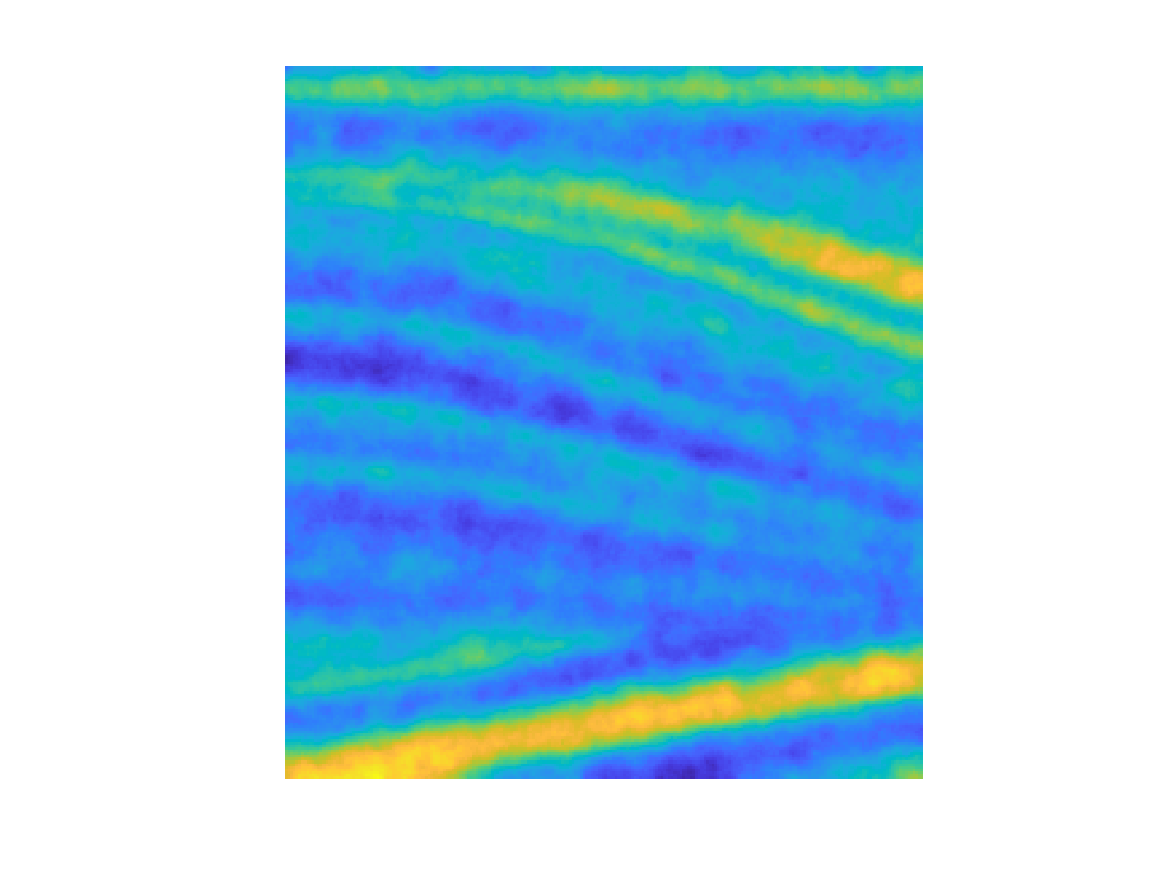} &
\hspace{-1.3cm}\includegraphics[width = 5cm]{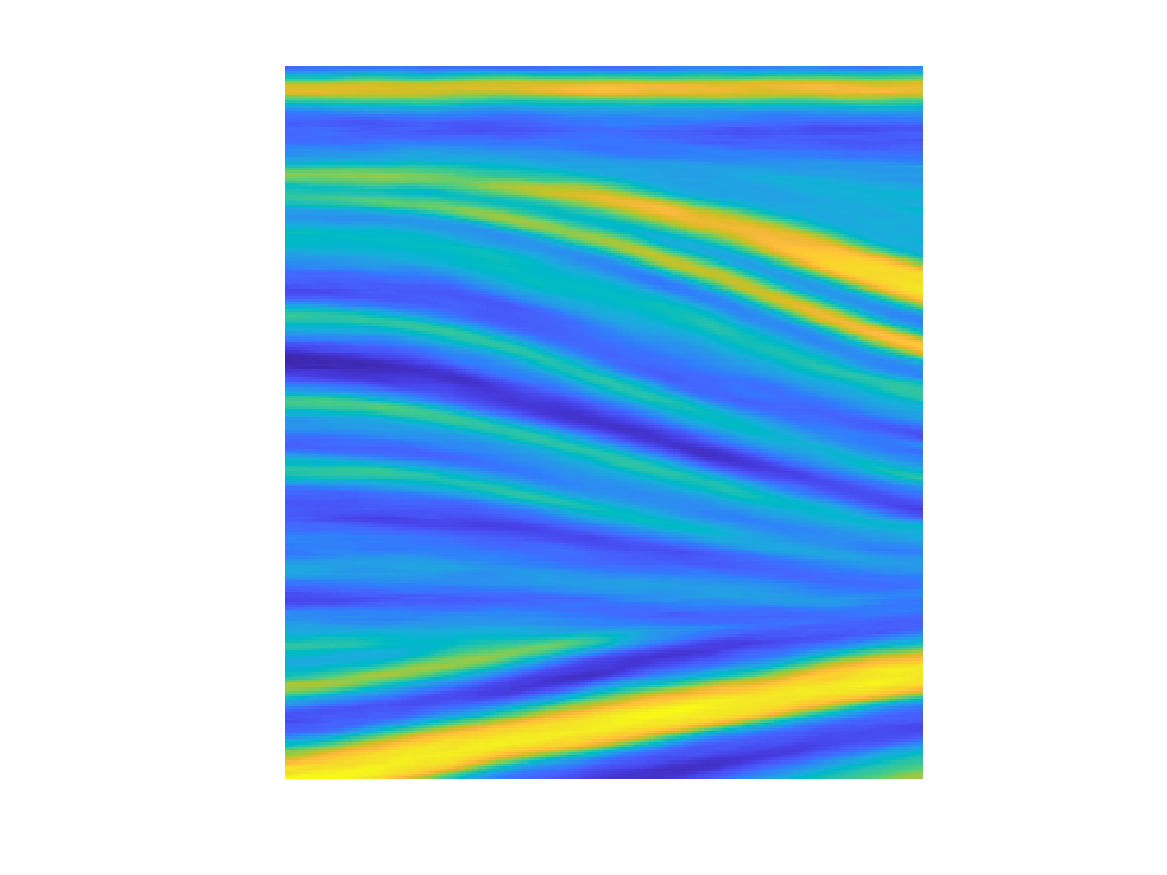}
\end{tabular}
\caption{Denoising test problem. (a) exact image; (b) available data; (c) image recovered by the (isotropic) Tikhonov regularization method; (d) image recovered by the new anisotropic Tikhonov regularization method.}
\label{fig: denoise_imgs}
\end{figure}
\end{center}

\begin{center}
\begin{figure}
\begin{tabular}{cc}
\hspace{-1.3cm}\small{(a)} & 
\hspace{0.3cm}\small{(b)}\\
\hspace{-1.3cm}\includegraphics[height = 5cm]{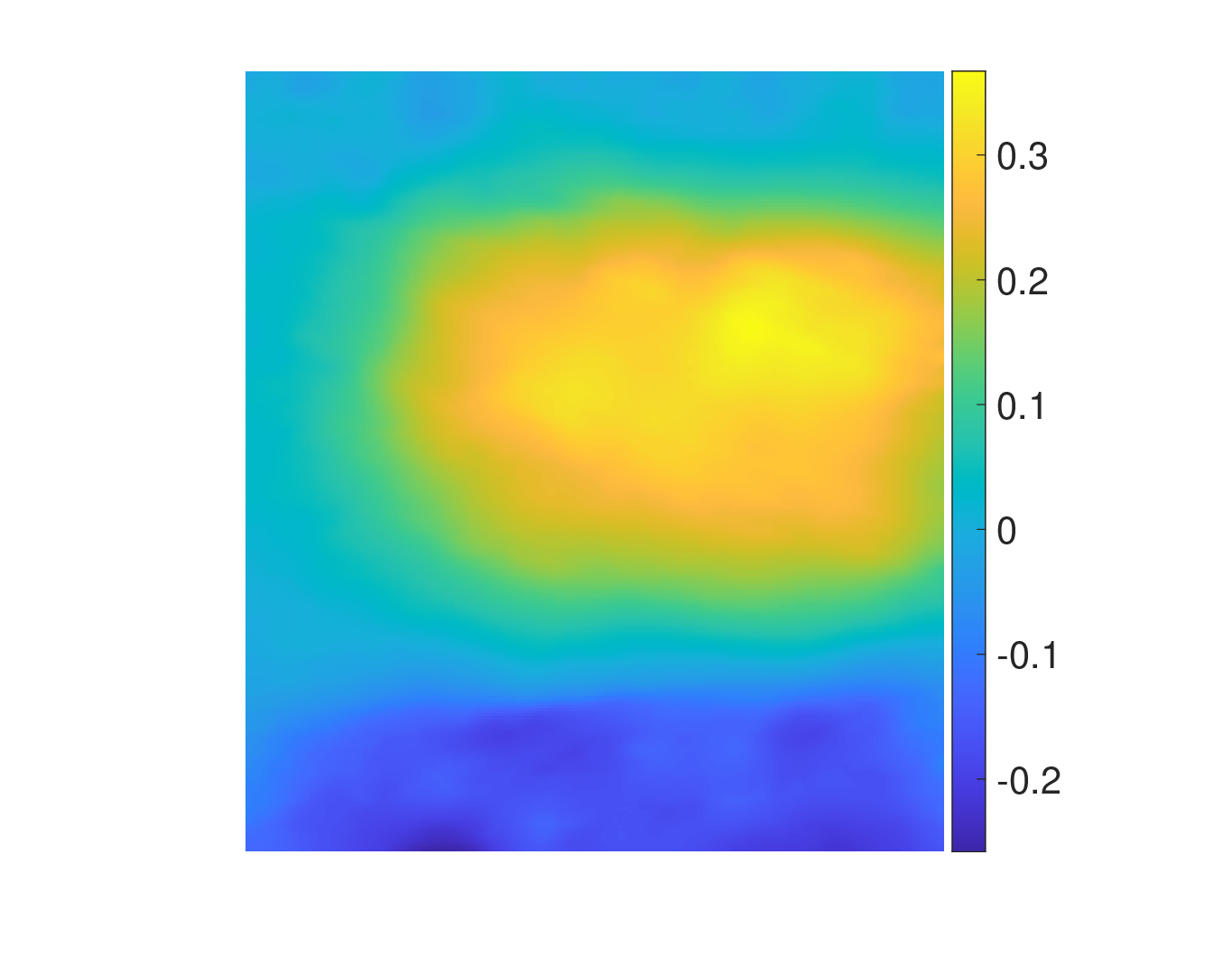} &
\hspace{0.3cm}\includegraphics[height = 5cm]{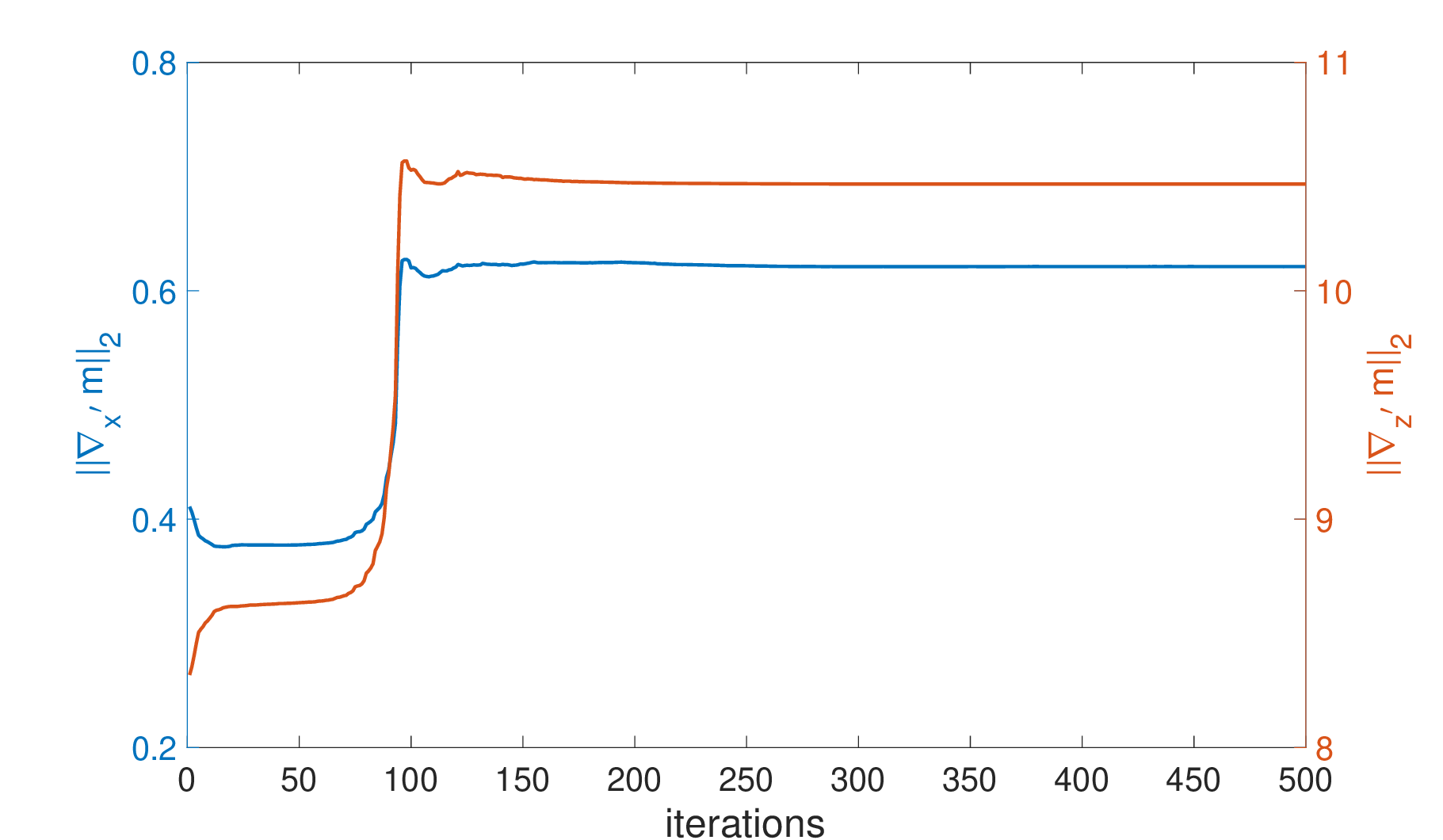}
\end{tabular}
\caption{Denoising test problem. (a) pixel-wise orientation parameters recovered solving the bilevel optimization method \eqref{eq:bilevel}; (b) 2-norm of the directional derivatives along $x'$ and $z'$ versus L-BFGS-B iterations.}
\label{fig: denoise_theta}
\end{figure}
\end{center}

\begin{center}
\begin{figure}
\begin{tabular}{ccc}
\hspace{-0.5cm}\small{(a)} & 
\hspace{-0.5cm}\small{(b)} &
\hspace{-0.5cm}\small{(c)}\\
\hspace{-0.5cm}\includegraphics[width = 5cm]{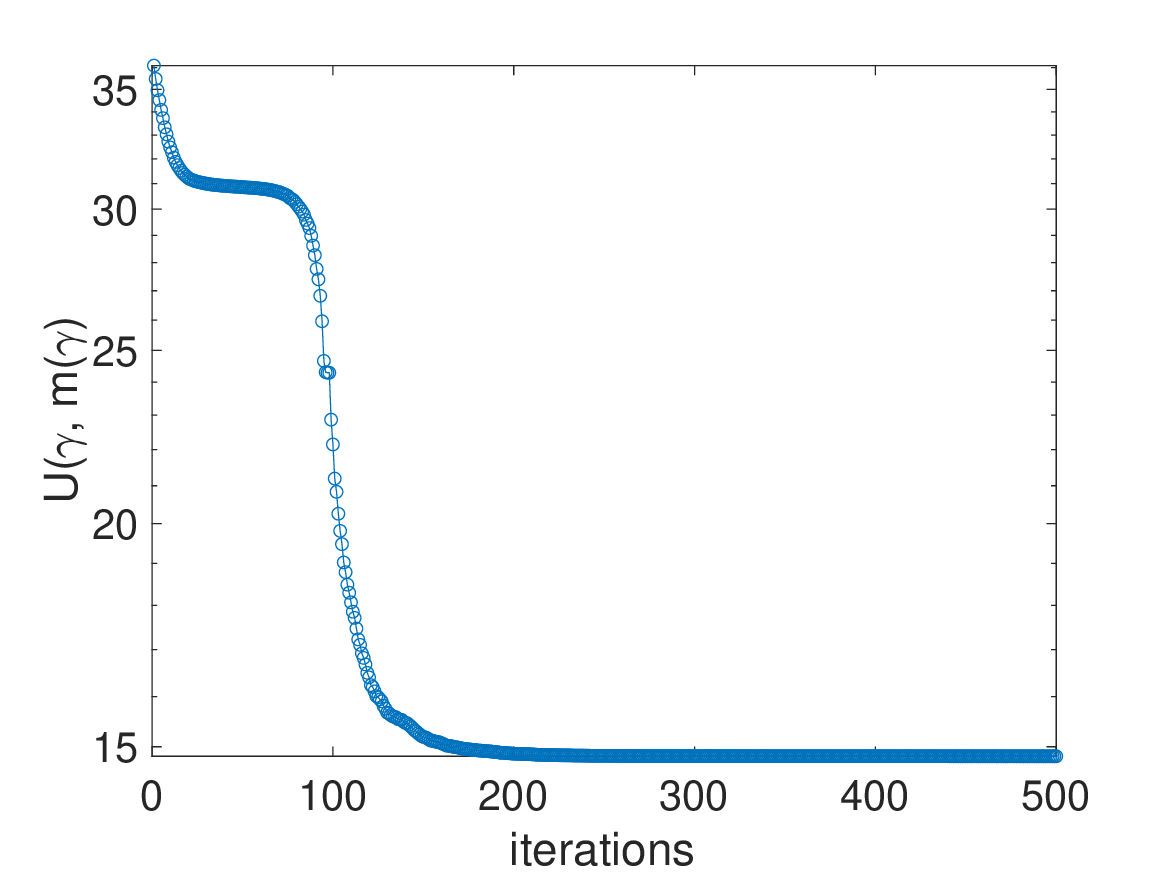} &
\hspace{-0.5cm}\includegraphics[width = 5cm]{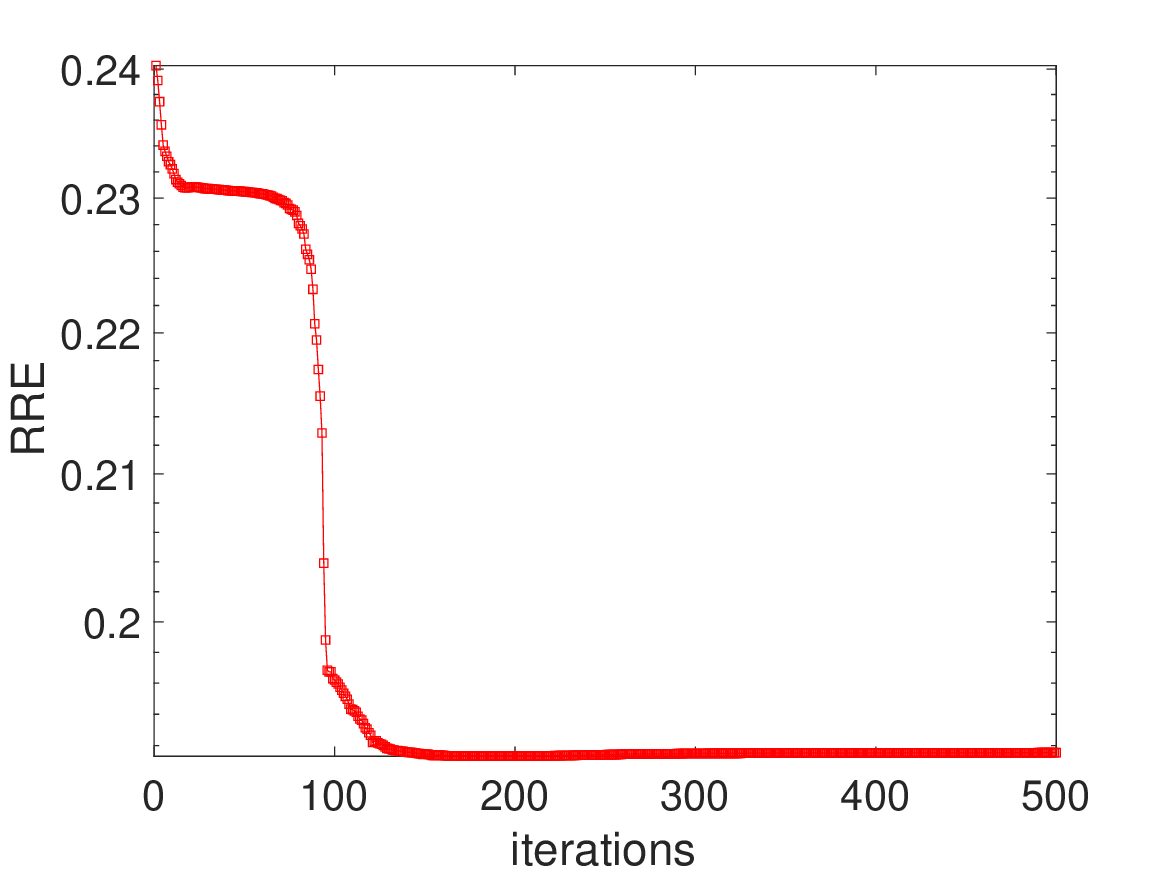} &
\hspace{-0.5cm}\includegraphics[width = 5cm]{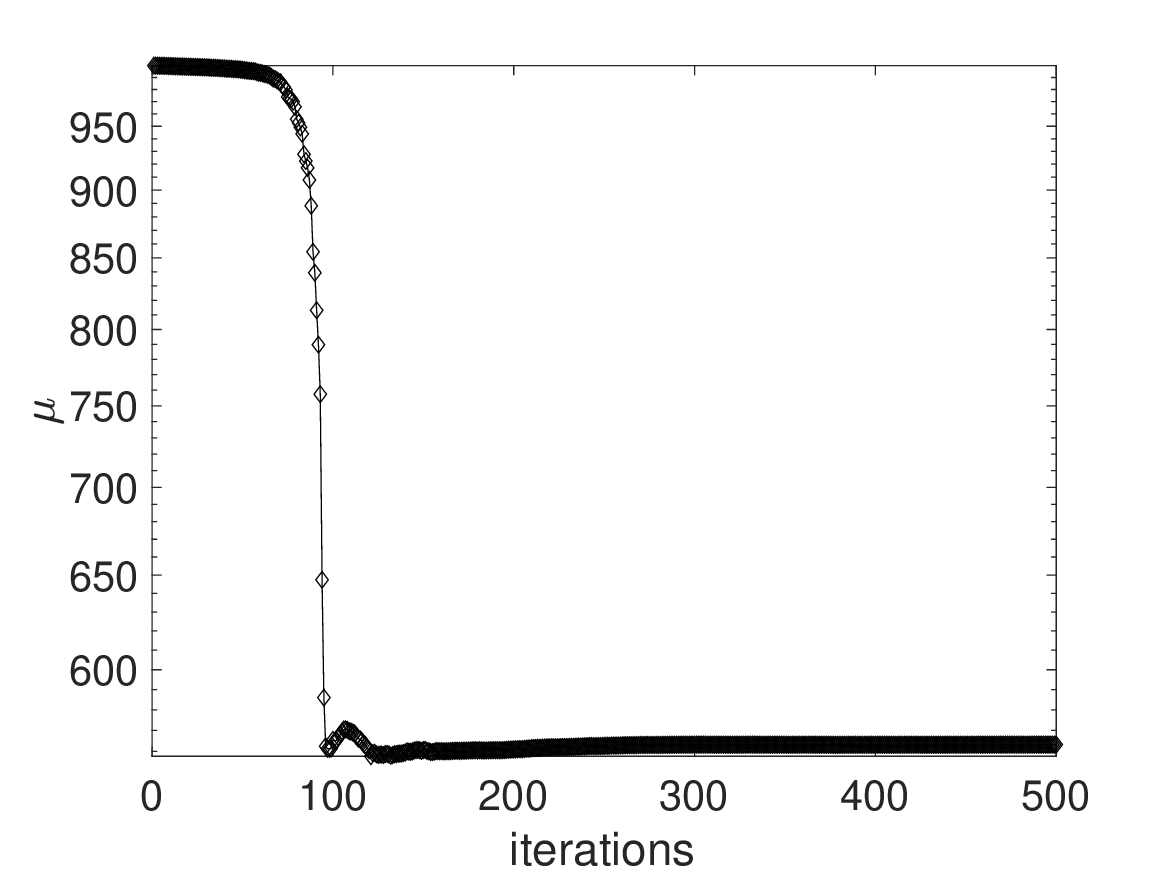}
\end{tabular}
\caption{Denoising test problem: (a) upper level objective function values versus L-BFGS-B iterations; (b) relative reconstruction errors versus L-BFGS-B iterations; (c) Tikhonov regularization parameter $\mu$ versus L-BFGS-B iterations.}
\label{fig: denoise_vals}
\end{figure}
\end{center}

We consider the restoration of the $238\times 266$ pixel image displayed in Figure \ref{fig: denoise_imgs}, frame (a). The corrupted image, displayed in Figure \ref{fig: denoise_imgs}, frame (b), is affected by some Gaussian noise $\be$ of level $\|\be\|_2/\|\bG\bm_{\rm true}\|_2$ equal to $1.6$. As weights for the gradient components we take $\sigma_1^{x'}=1$ and $\sigma_2^{z'} = 10^{-1}$. The parameters appearing in the upper level objective function \eqref{ULobjfcn} are taken as $\alpha=10$, $\beta=15$. The reconstructions by (isotropic) Tikhonov regulatization and the new local anisotropic regularization strategy are displayed in Figure \ref{fig: denoise_imgs}, frames (c) and (d), respectively. The regularization parameter for (isotropic) Tikhonov regularization recovered by the discrepancy principle is $\mu=26.21$. The relative reconstruction error associated to such method is $0.3171$. The local orientation parameter recovered by the new bilevel optimization method are displayed in Figure \ref{fig: denoise_theta}. The history of relevant quantities that monitor the progress of the new bilevel optimization approach to bilevel optimization are displayed in Figure \ref{fig: denoise_vals}.

\subsection{Image deblurring}
\begin{center}
\begin{figure}
\begin{tabular}{cccc}
\hspace{-1.3cm}\small{(a)} & 
\hspace{-1.3cm}\small{(b)} & 
\hspace{-1.3cm}\small{(c)} & 
\hspace{-1.3cm}\small{(d)} \vspace{-0.1cm}\\
\hspace{-1.3cm}\includegraphics[width = 5cm]{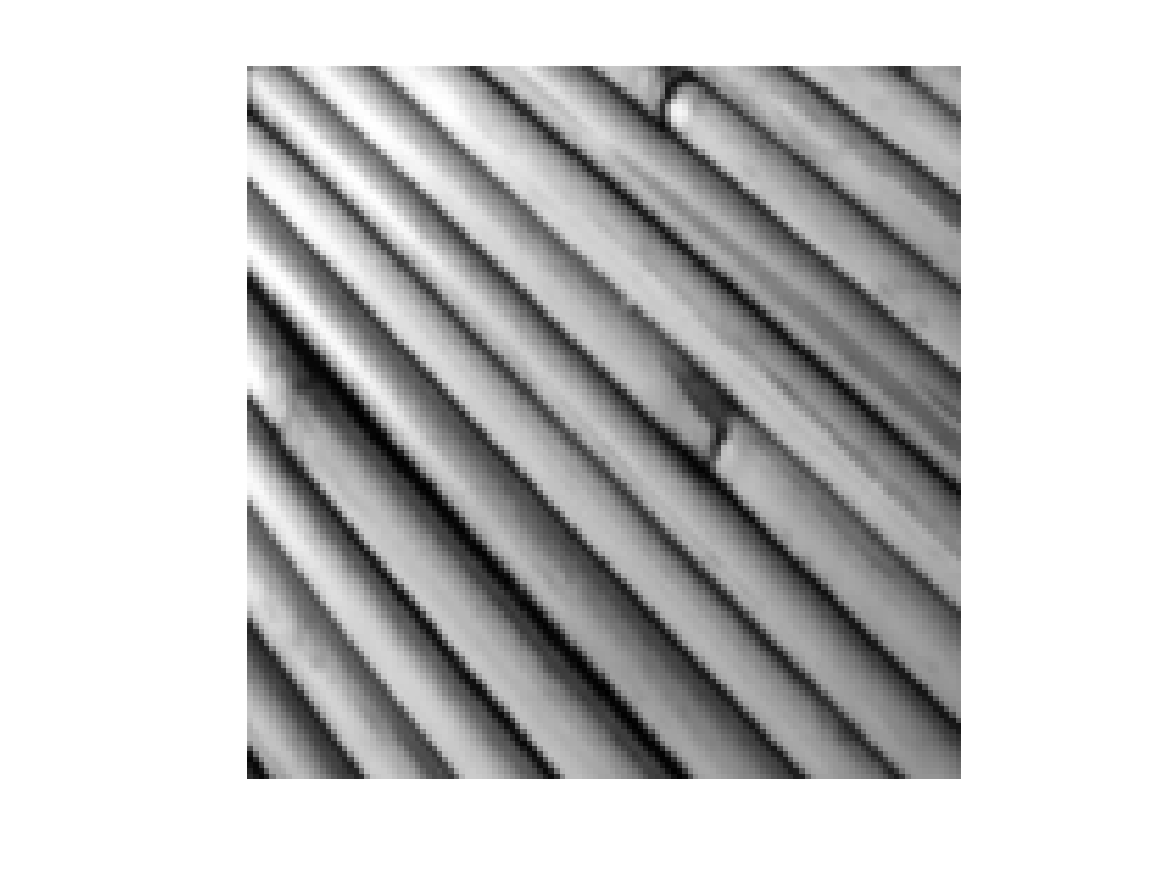} &
\hspace{-1.3cm}\includegraphics[width = 5cm]{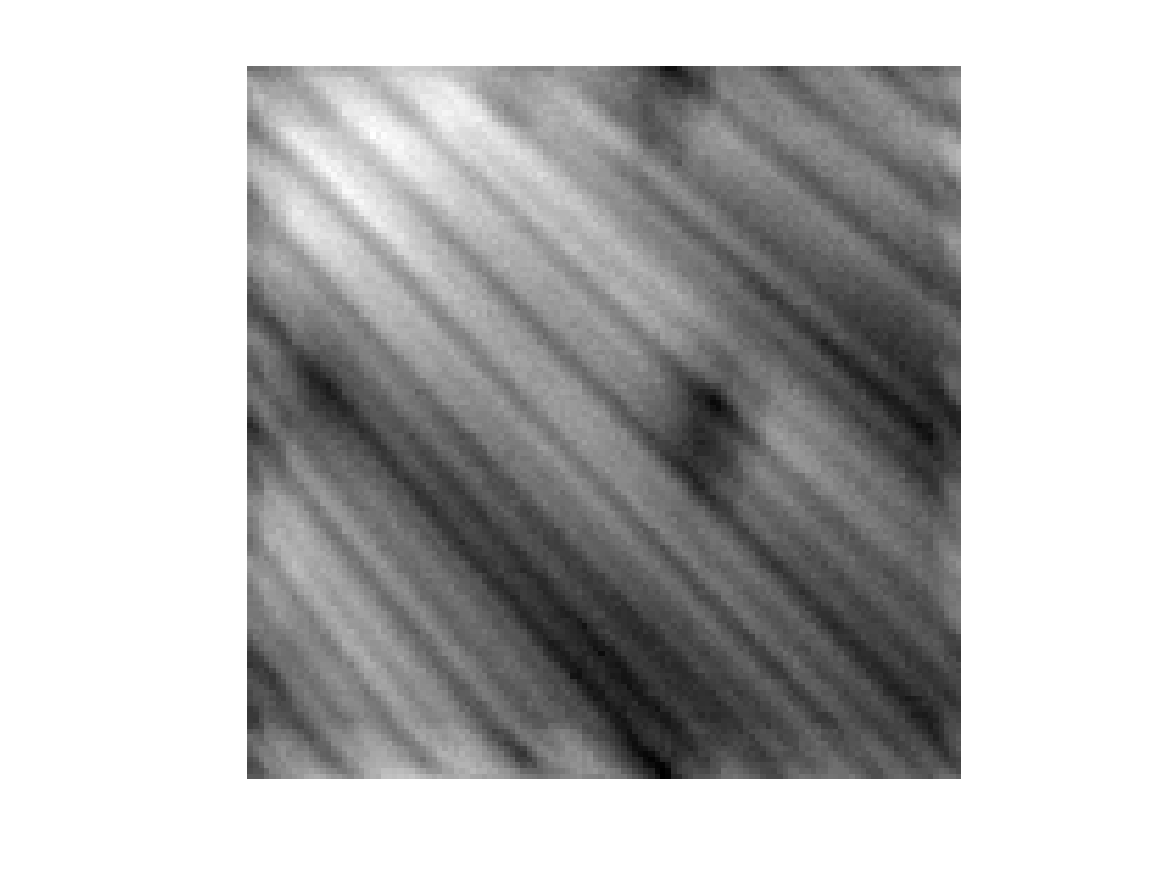} & 
\hspace{-1.3cm}\includegraphics[width = 5cm]{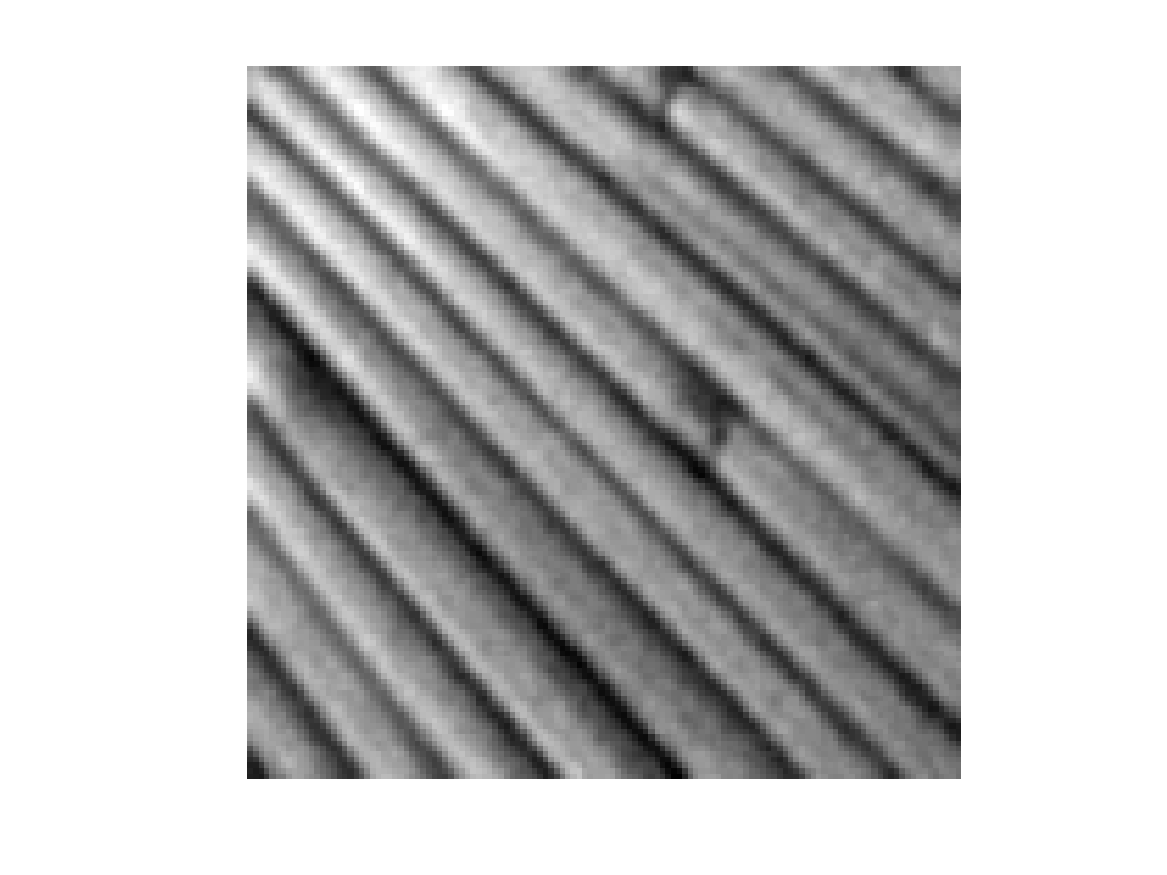} &
\hspace{-1.3cm}\includegraphics[width = 5cm]{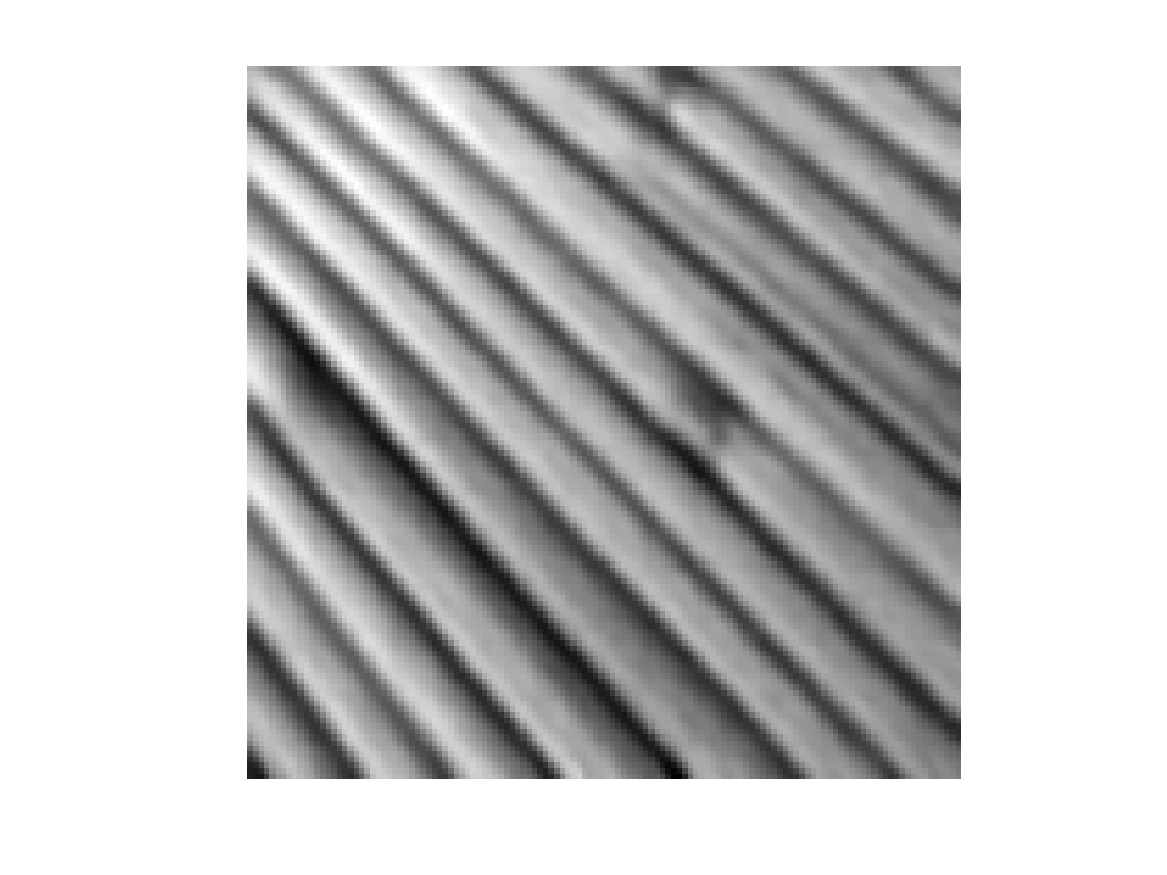}
\end{tabular}
\caption{Deblurring test problem. (a) exact image; (b) available data; (c) image recovered by the (isotropic) Tikhonov regularization method; (d) image recovered by the new anisotropic Tikhonov regularization method.}
\label{fig: deblur_imgs}
\end{figure}
\end{center}

\begin{center}
\begin{figure}
\begin{tabular}{cc}
\hspace{-1.3cm}\small{(a)} & 
\hspace{0.3cm}\small{(b)}\\
\hspace{-1.3cm}\includegraphics[height = 5cm]{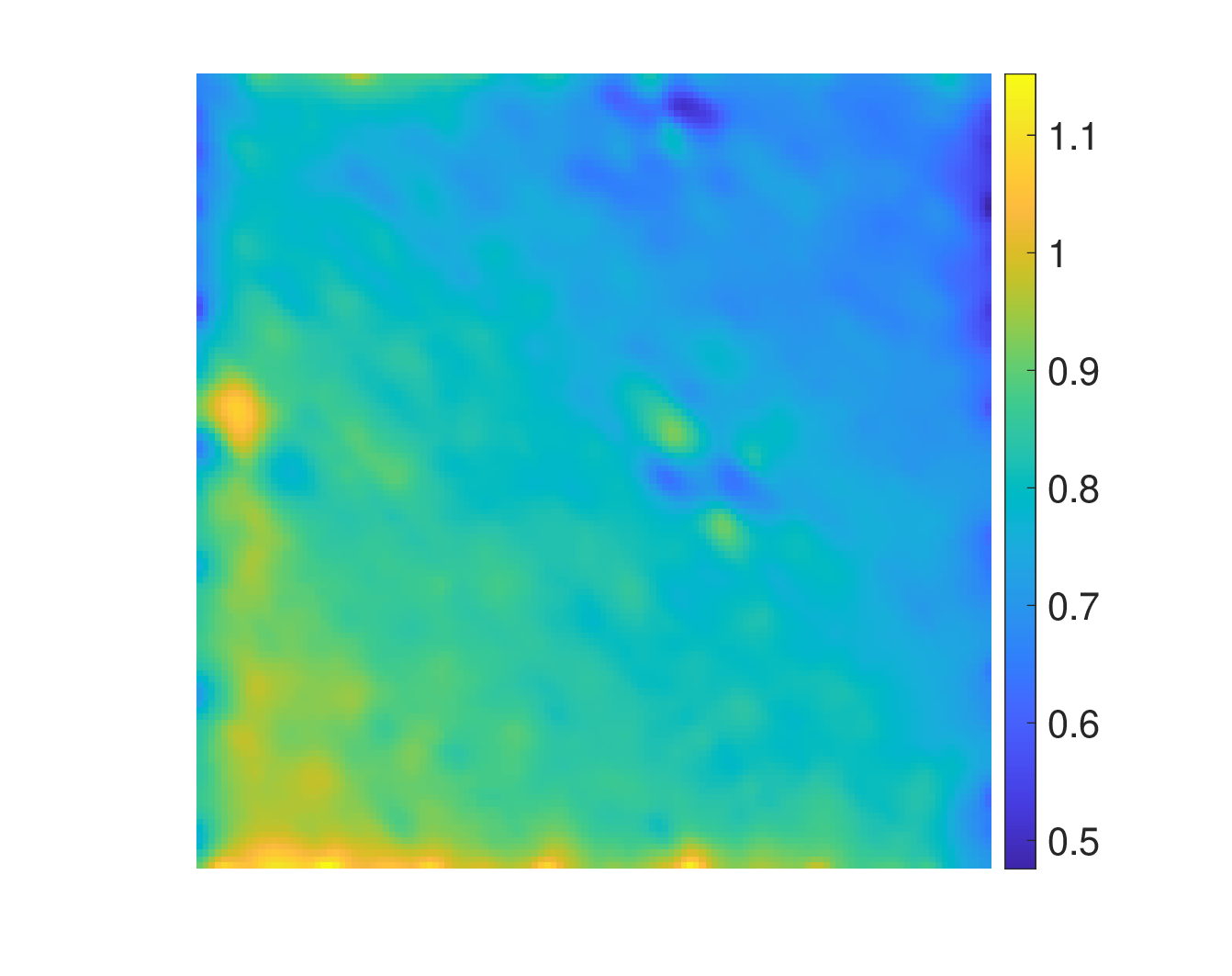} &
\hspace{0.3cm}\includegraphics[height = 5cm]{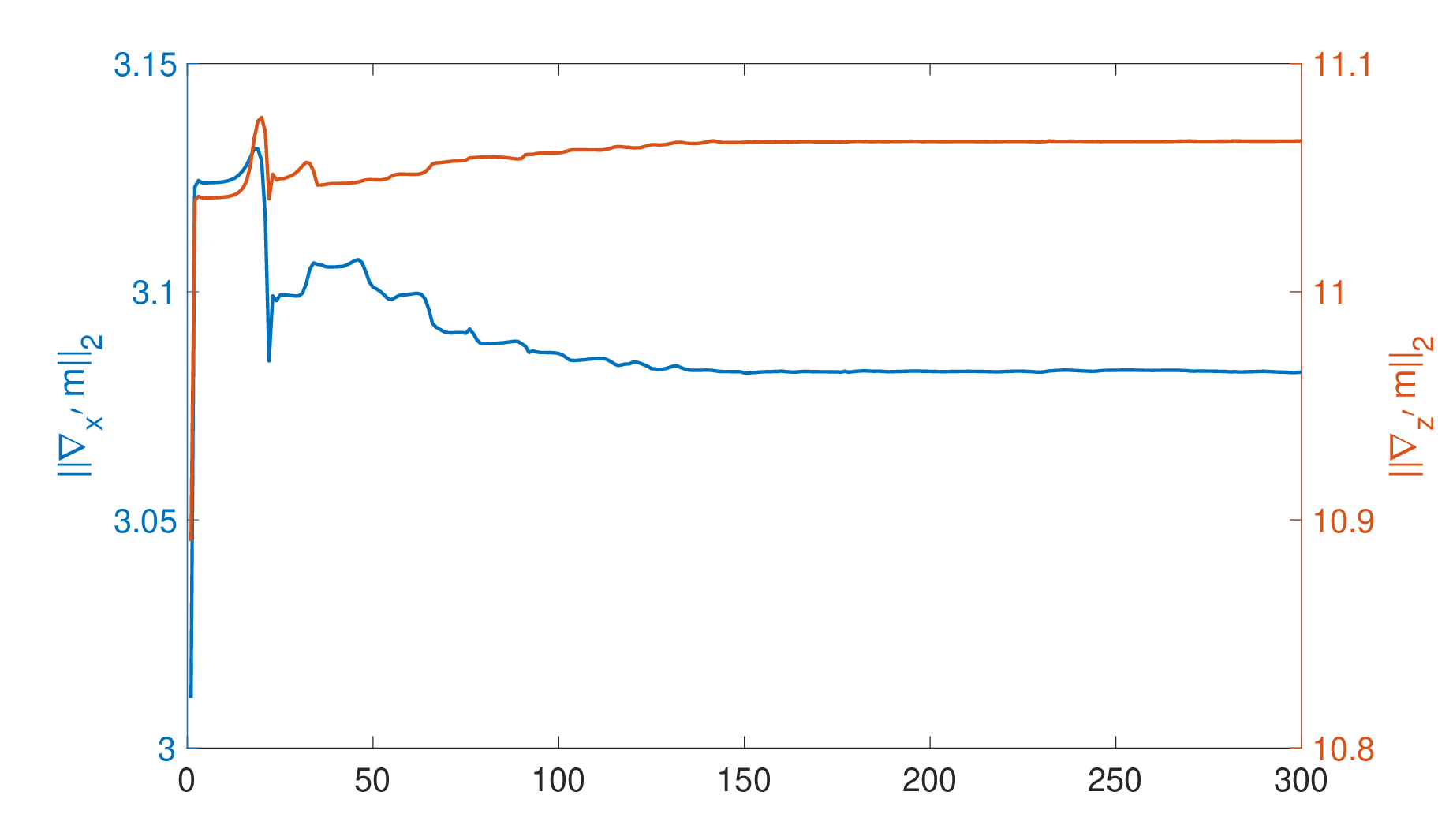}
\end{tabular}
\caption{Deblurring test problem. (a) pixel-wise orientation parameters recovered solving the bilevel optimization method \eqref{eq:bilevel}; (b) 2-norm of the directional derivatives along $x'$ and $z'$ versus L-BFGS-B iterations.}
\label{fig: deblur_theta}
\end{figure}
\end{center}

\begin{center}
\begin{figure}
\begin{tabular}{ccc}
\hspace{-0.5cm}\small{(a)} & 
\hspace{-0.5cm}\small{(b)} &
\hspace{-0.5cm}\small{(c)}\\
\hspace{-0.5cm}\includegraphics[width = 5cm]{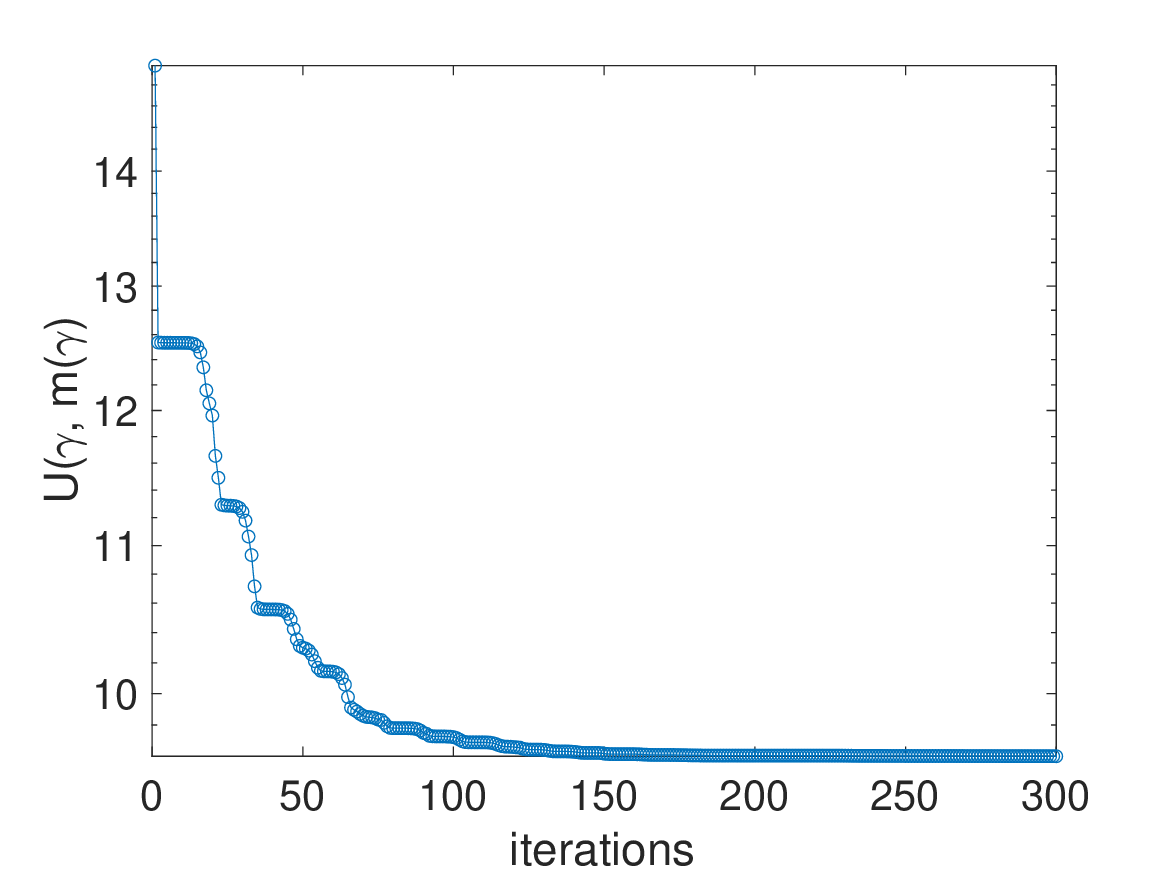} &
\hspace{-0.5cm}\includegraphics[width = 5cm]{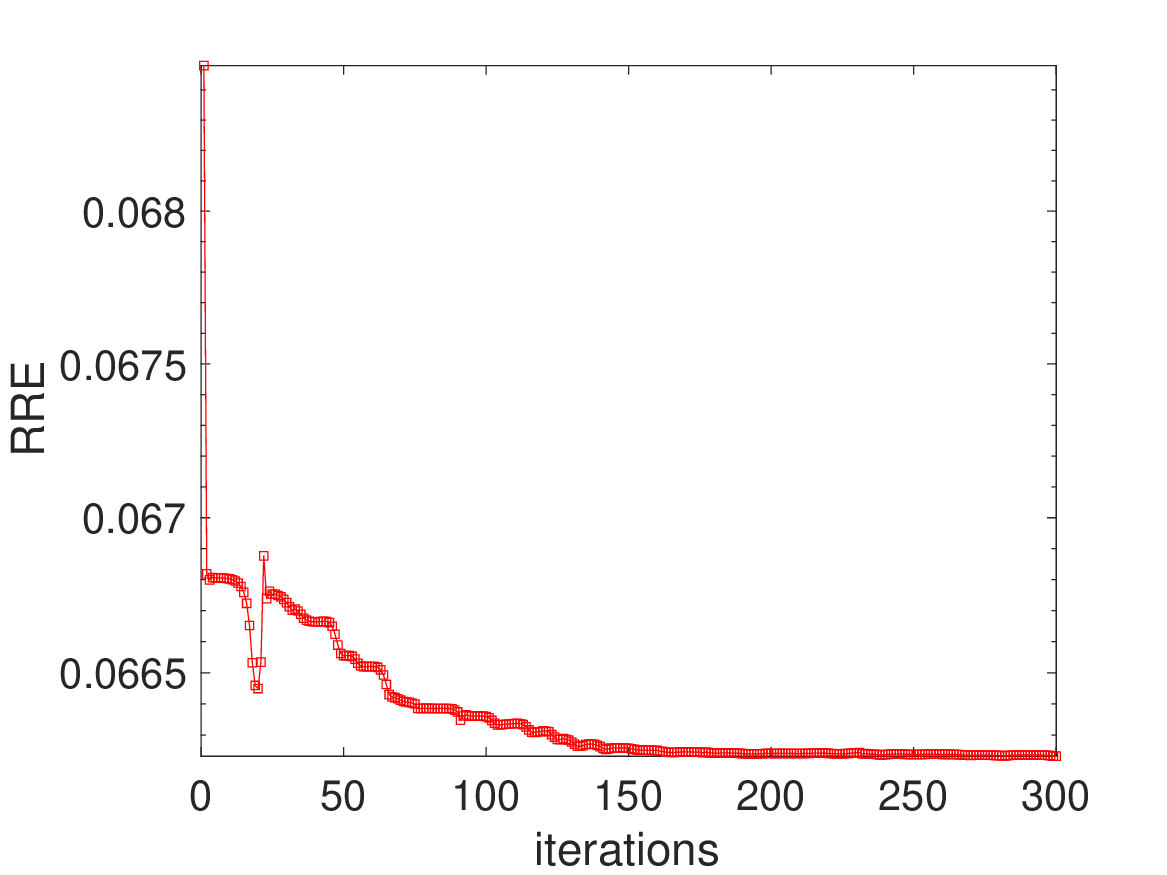} &
\hspace{-0.5cm}\includegraphics[width = 5cm]{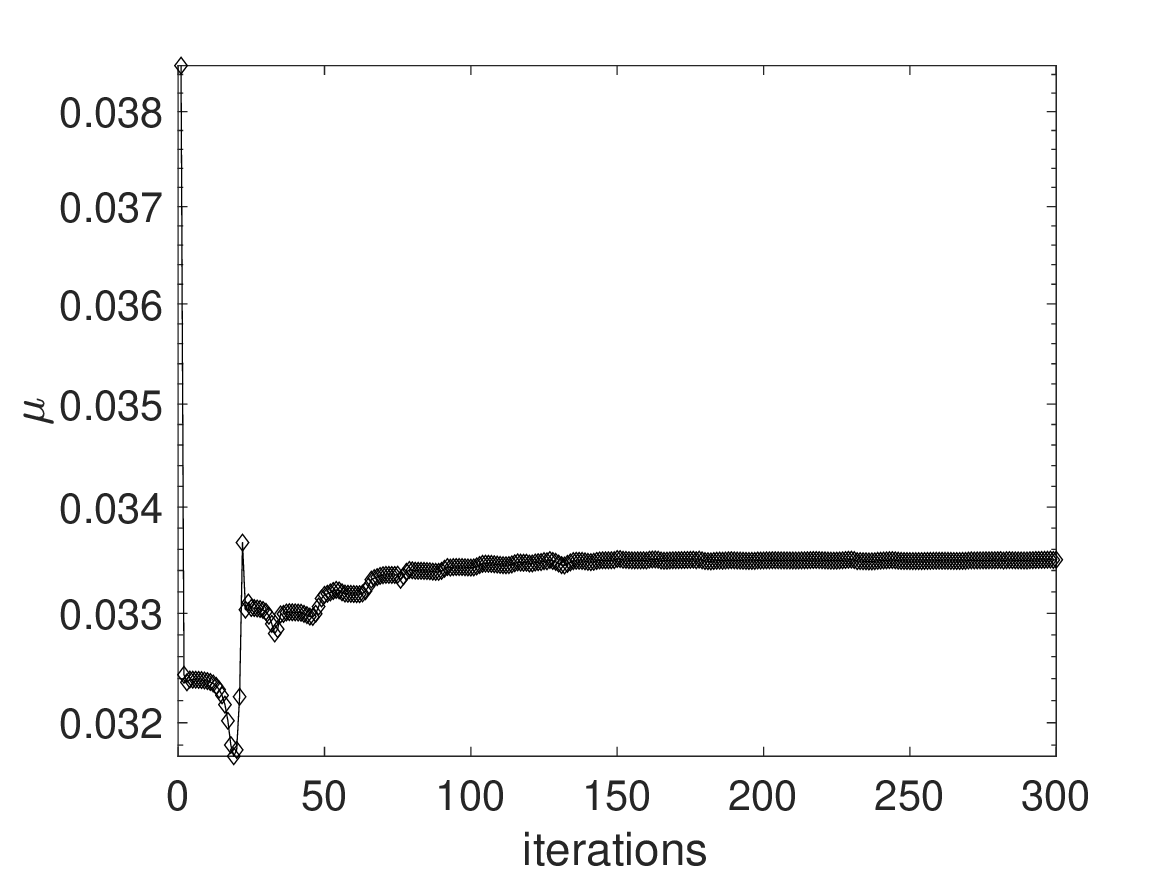}
\end{tabular}
\caption{Deblurring test problem. (a) upper level objective function values versus L-BFGS-B iterations; (b) relative reconstruction errors versus L-BFGS-B iterations; (c) Tikhonov regularization parameter $\mu$ versus L-BFGS-B iterations.}
\label{fig: deblur_vals}
\end{figure}
\end{center}

We consider the restoration of a $128\times 128$ pixel image of a bamboo fence, which has undergone a blurring process, with a Gaussian PSF whose standard deviation is 36; some Gaussian noise $\be$ of level $\|\be\|_2/\|\bG\bm_{\rm true}\|_2$ equal to $1\cdot 10^{-2}$ has been added to the blurred image. As weights for the gradient components we take $\sigma_1^{x'}=1$ and $\sigma_2^{z'} = 10^{-3}$. The parameters appearing in the upper level objective function \eqref{ULobjfcn} are taken as $\alpha=\beta=4\cdot 10^{-3}$. The considered ground truth and data are displayed in Figure \ref{fig: deblur_imgs}, along with the reconstructions by (isotropic) Tikhonov regulatization and the new local anisotropic regularization strategy. The regularization parameter for (isotropic) Tikhonov regularization recovered by the discrepancy principle is $\mu=0.0061$. The relative reconstruction error associated to the former is $0.0882$. The local orientation parameter recovered by the new bilevel optimization method are displayed in Figure \ref{fig: deblur_theta}. The history of relevant quantities that monitor the progress of the new bilevel optimization approach to bilevel optimization are displayed in Figure \ref{fig: deblur_vals}.

\subsection{Seismic tomography} 

\begin{center}
\begin{figure}
\begin{tabular}{cccc}
\hspace{-1.3cm}\small{(a)} & 
\hspace{-1.3cm}\small{(b)} & 
\hspace{-1.3cm}\small{(c)} & 
\hspace{-1.3cm}\small{(d)} \vspace{-0.1cm}\\
\hspace{-1.3cm}\includegraphics[width = 5cm]{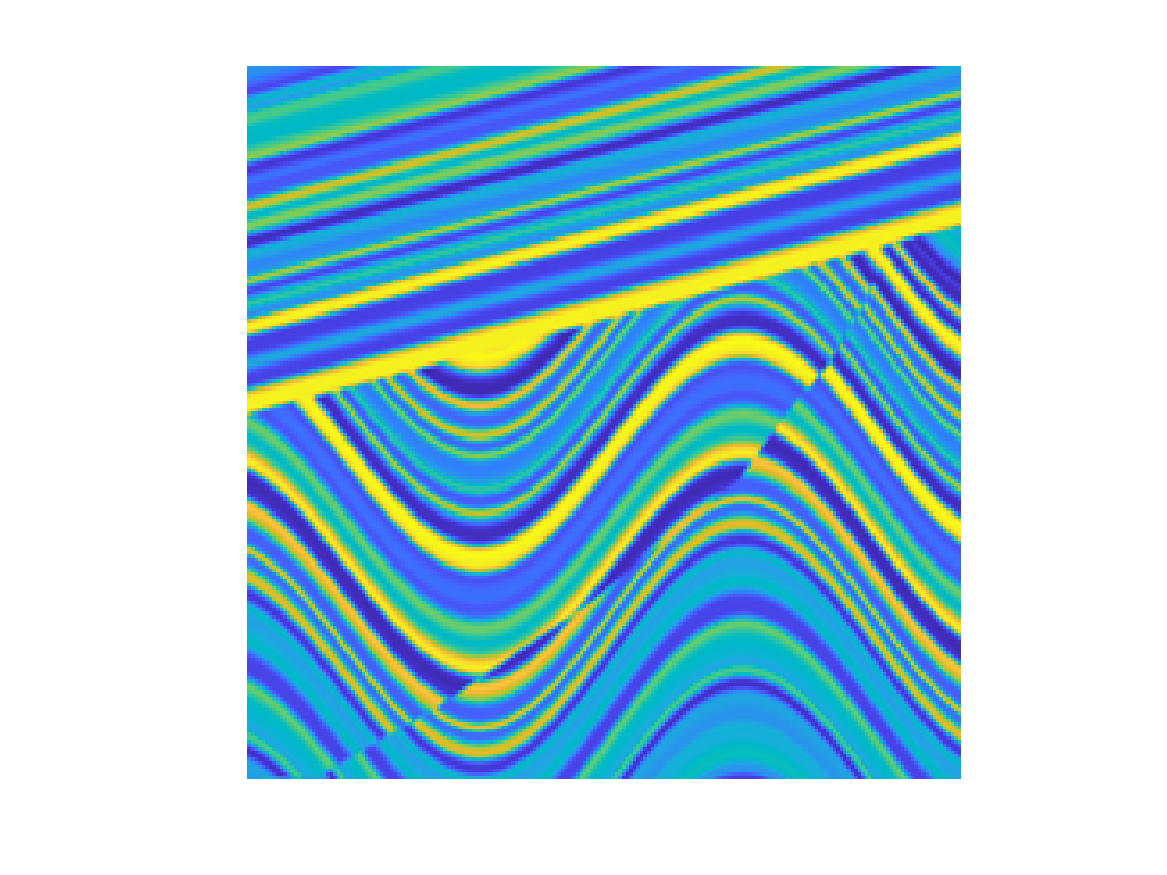} &
\hspace{-1.3cm}\includegraphics[width = 5cm]{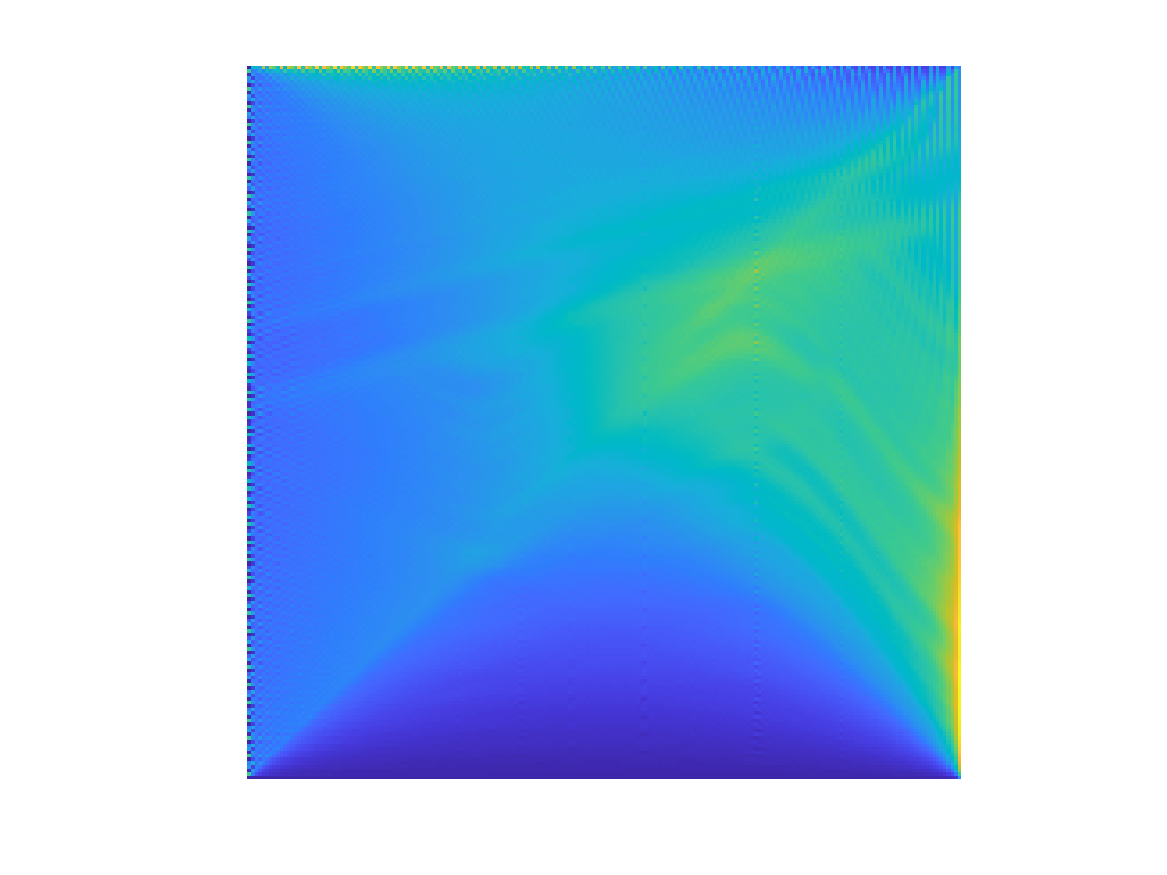} & 
\hspace{-1.3cm}\includegraphics[width = 5cm]{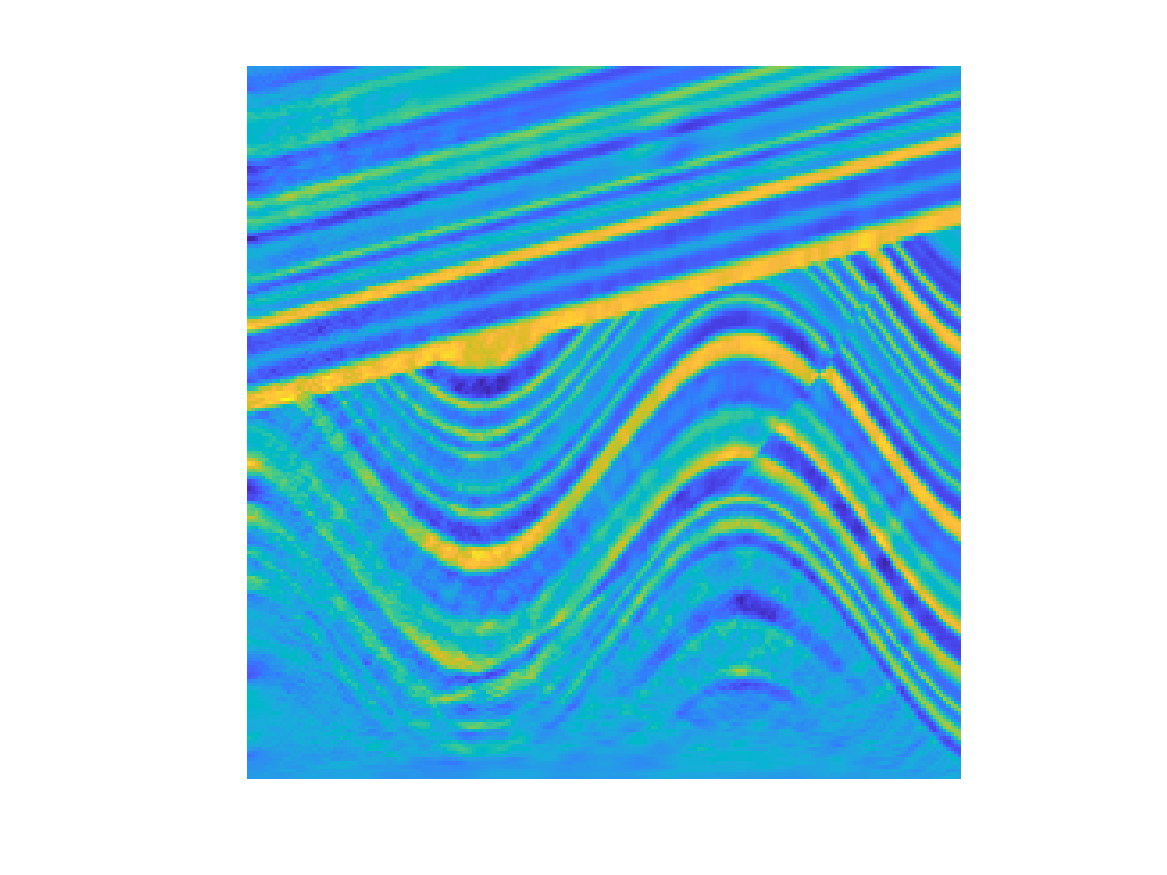} &
\hspace{-1.3cm}\includegraphics[width = 5cm]{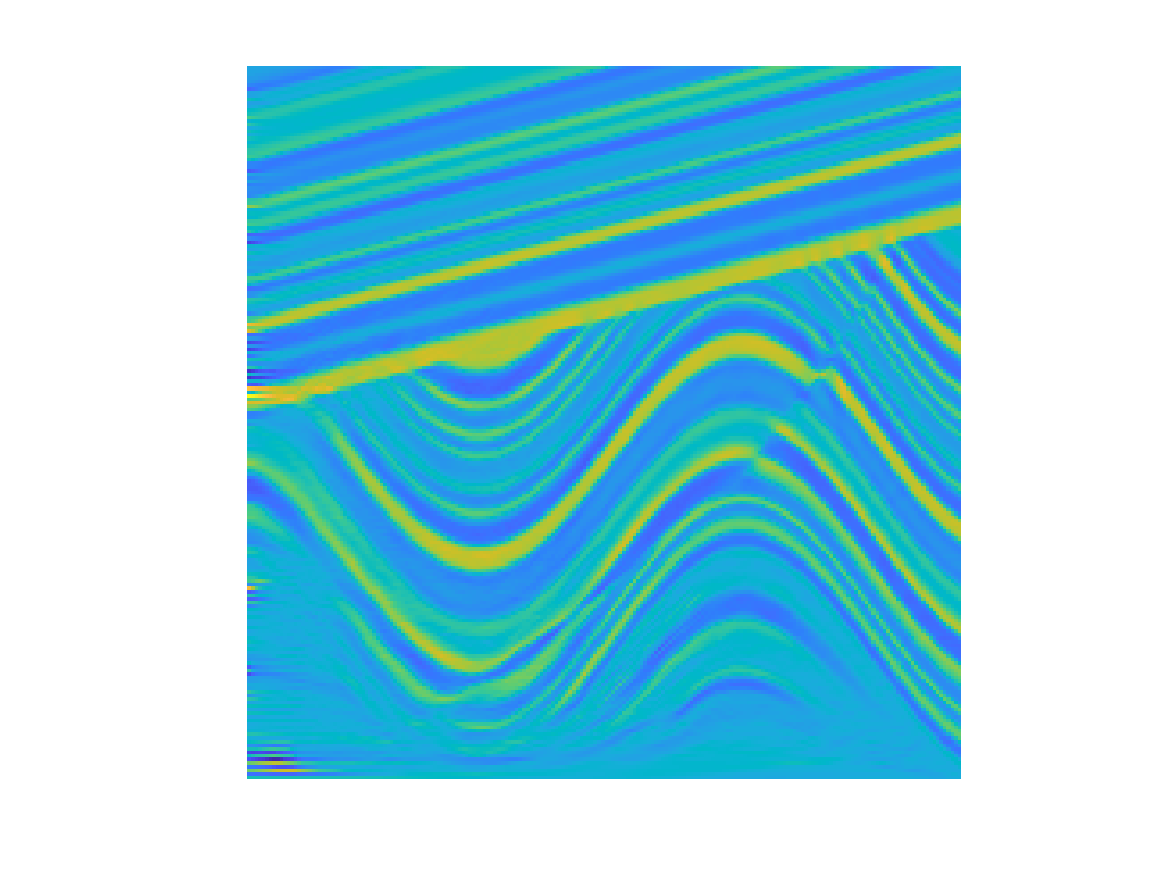}
\end{tabular}
\caption{Seismic tomography test problem. (a) exact phantom; (b) phantom recovered by backprojection; (c) phantom recovered by the (isotropic) Tikhonov regularization method; (d) phantom recovered by the new anisotropic Tikhonov regularization method.}
\label{fig: tomo_imgs}
\end{figure}
\end{center}

\begin{center}
\begin{figure}
\begin{tabular}{cc}
\hspace{-1.3cm}\small{(a)} & 
\hspace{0.3cm}\small{(b)}\\
\hspace{-1.3cm}\includegraphics[height = 5cm]{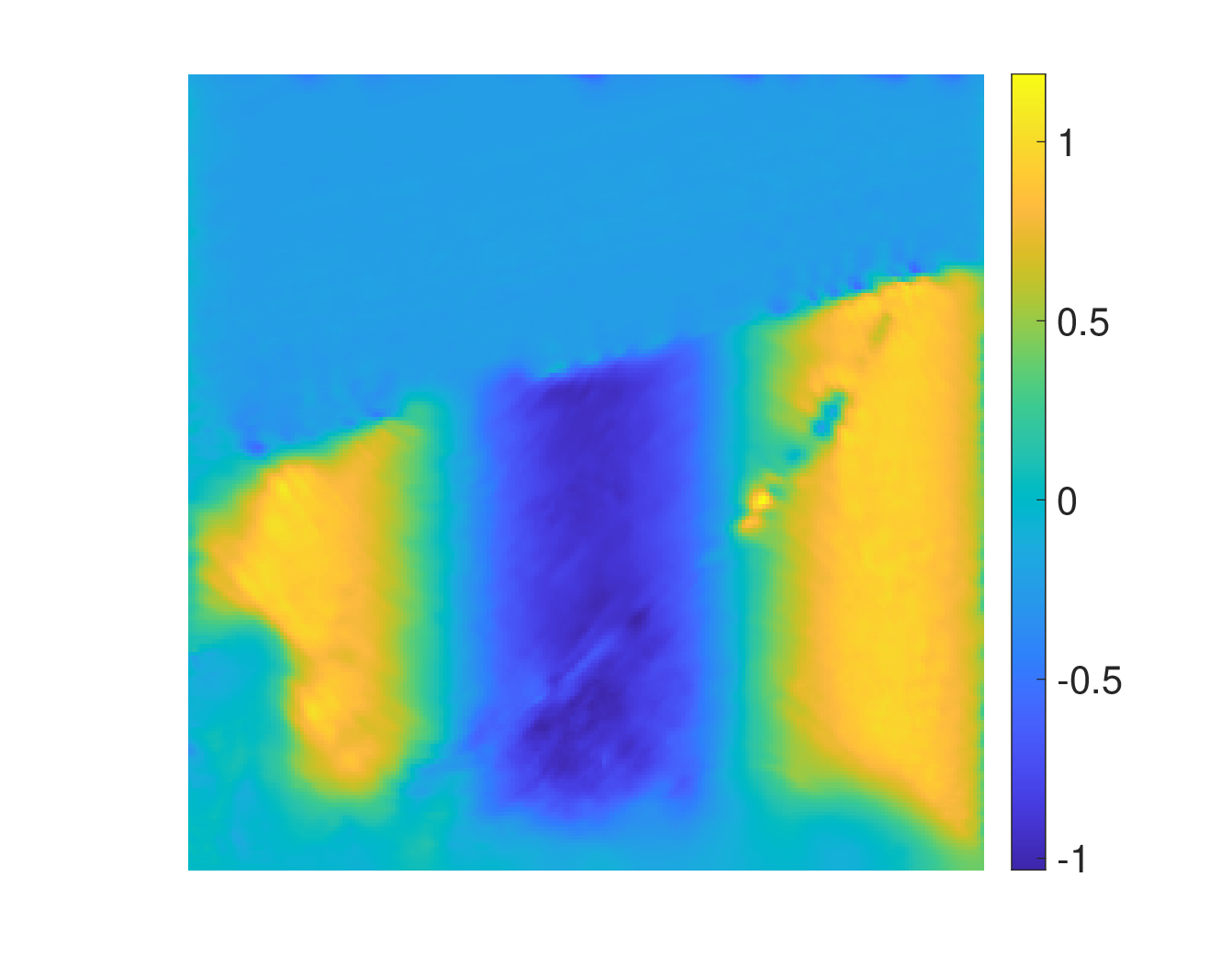} &
\hspace{0.3cm}\includegraphics[height = 5cm]{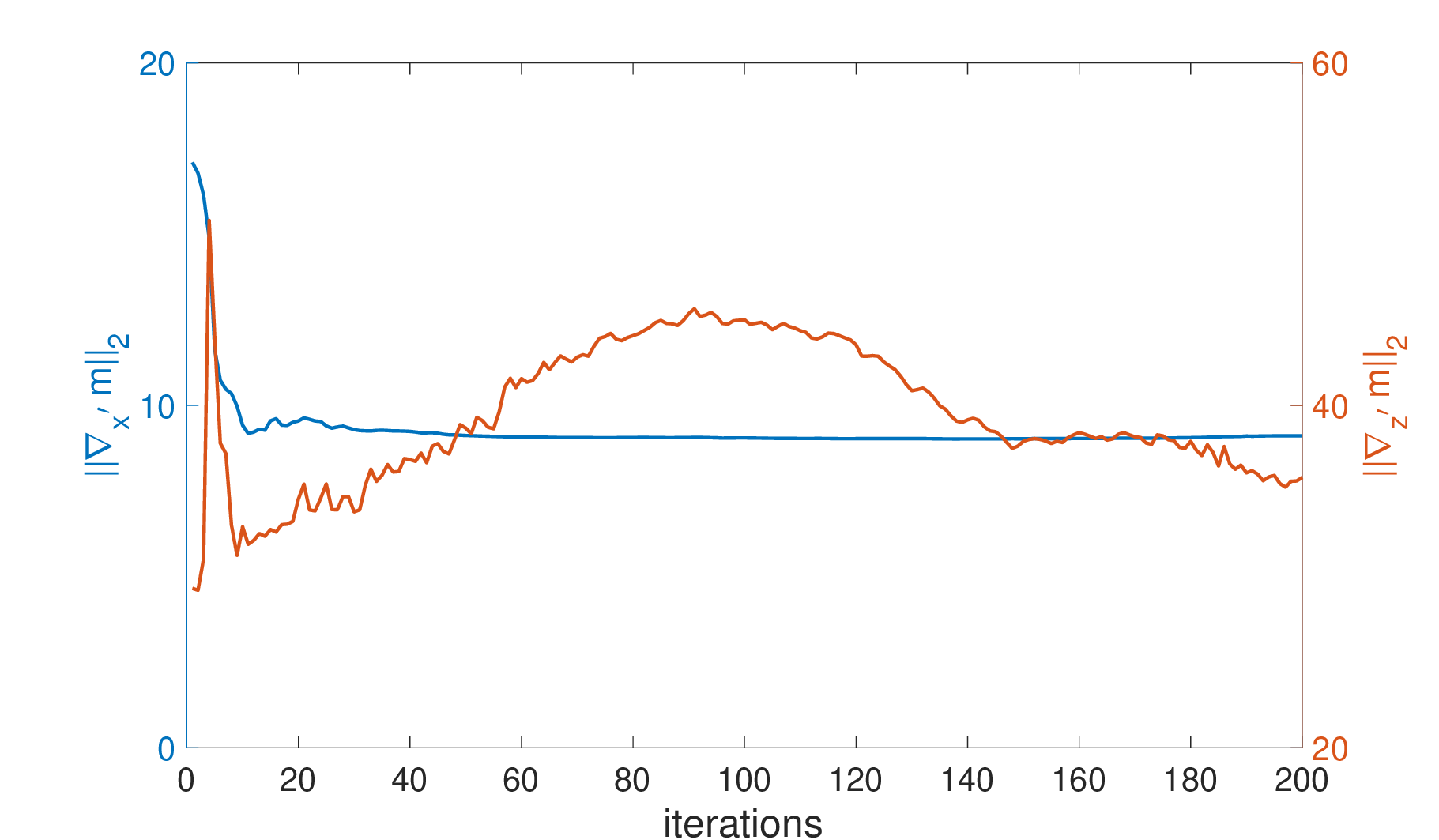}
\end{tabular}
\caption{Seismic tomography test problem. (a) pixel-wise orientation parameters recovered solving the bilevel optimization method \eqref{eq:bilevel}; (b) 2-norm of the directional derivatives along $x'$ and $z'$ versus L-BFGS-B iterations.}
\label{fig: tomo_theta}
\end{figure}
\end{center}

\begin{center}
\begin{figure}
\begin{tabular}{ccc}
\hspace{-0.5cm}\small{(a)} & 
\hspace{-0.5cm}\small{(b)} &
\hspace{-0.5cm}\small{(c)}\\
\hspace{-0.5cm}\includegraphics[width = 5cm]{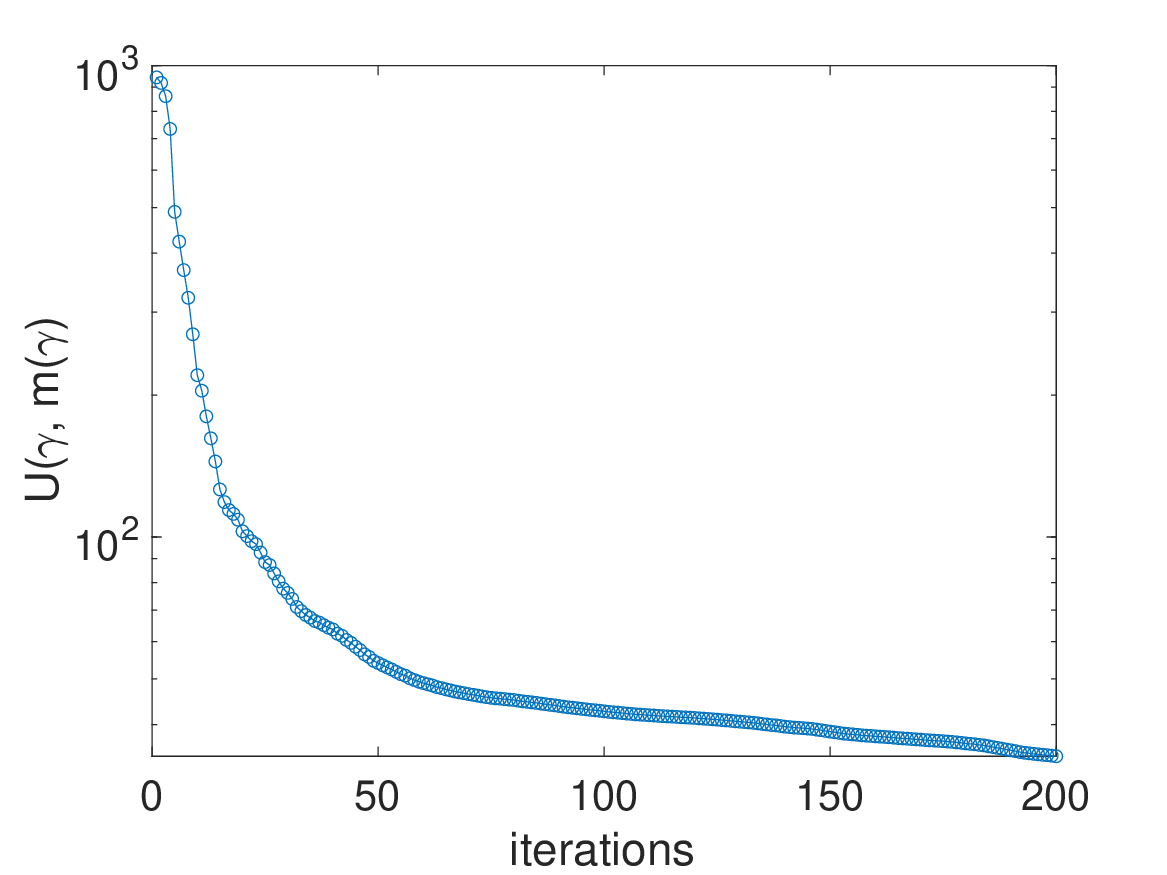} &
\hspace{-0.5cm}\includegraphics[width = 5cm]{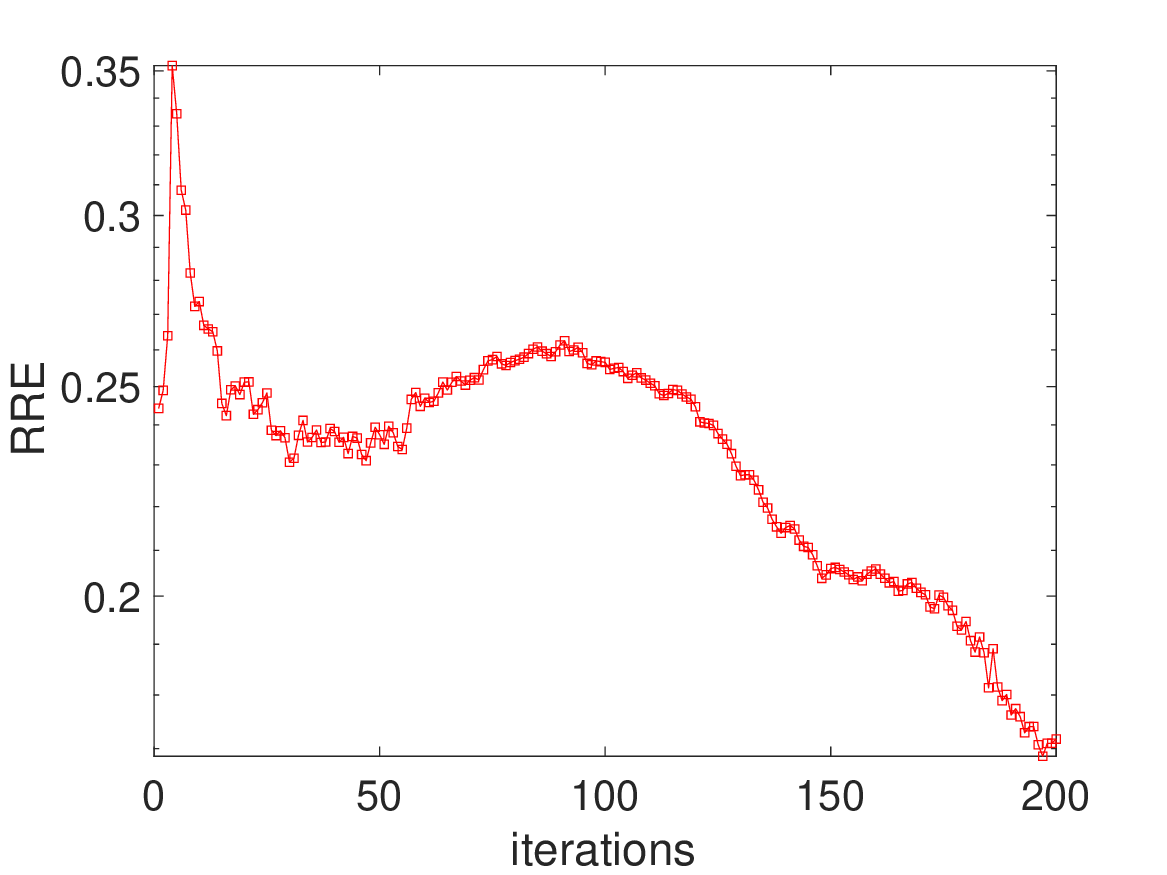} &
\hspace{-0.5cm}\includegraphics[width = 5cm]{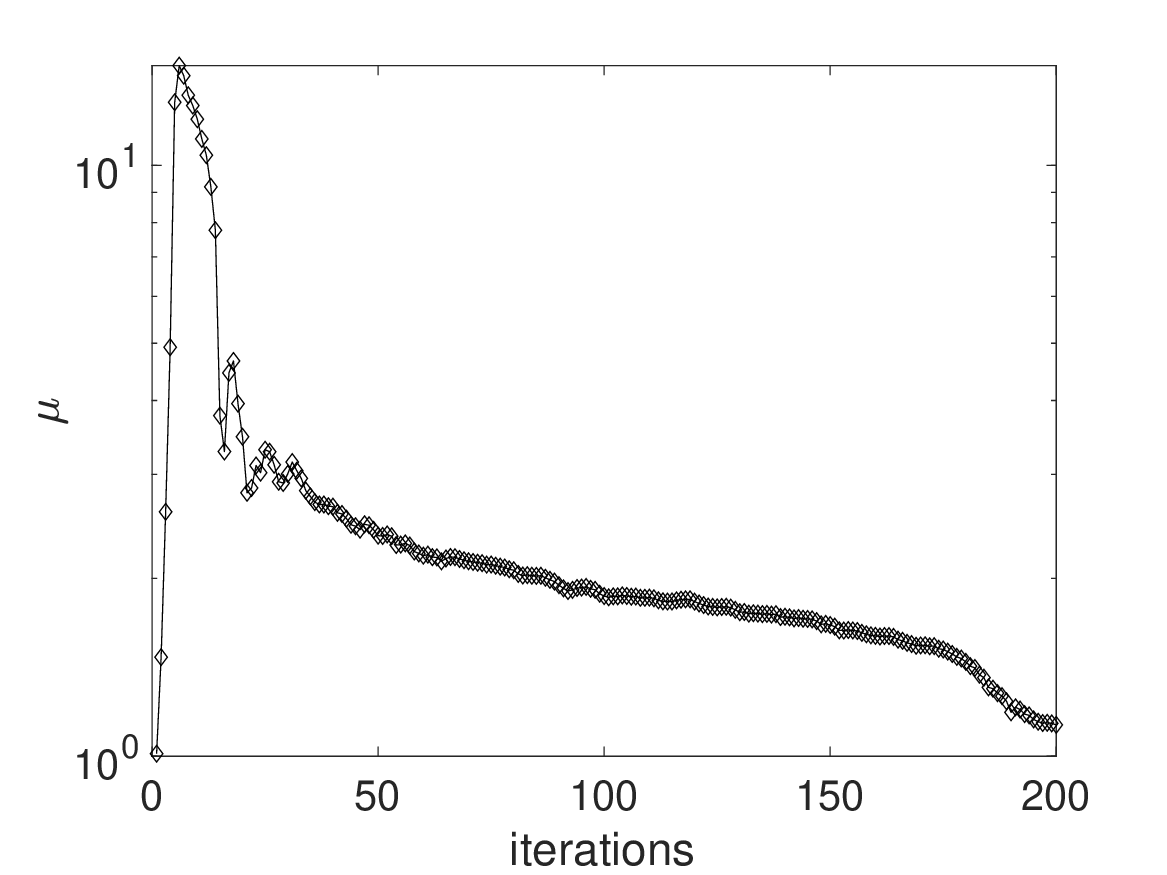}
\end{tabular}
\caption{Seismic tomography test problem. (a) upper level objective function values versus L-BFGS-B iterations; (b) relative reconstruction errors versus L-BFGS-B iterations; (c) Tikhonov regularization parameter $\mu$ versus L-BFGS-B iterations.}
\label{fig: tomo_vals}
\end{figure}
\end{center}

This test problem models the recovering of attenuation coefficients within a bounded domain, where rays emitted by 100 sources located at the right-hand side of the domain are recorded by 160 receivers (seismographs) equally spaced along the left and top
boundaries of the domain. The model parameter we wish to approximate are represented as a $200\times 200$ pixel image. The acquired data are affected by some Gaussian noise $\be$ of level $\|\be\|_2/\|\bG\bm_{\rm true}\|_2$ equal to $2.5\cdot 10^{-5}$. As weights for the gradient components we take $\sigma^{x'}=1$ and $\sigma^{z'} = 10^{-3}$. The parameters appearing in the upper level objective function \eqref{ULobjfcn} are taken as $\alpha=1$, $\beta=3\cdot 10^{-1}$. The considered ground truth, along with the reconstructions by simple backprojection, by (isotropic) Tikhonov regulatization and by the new local anisotropic regularization strategy, are displayed in Figure \ref{fig: tomo_imgs}. The regularization parameter for (isotropic) Tikhonov regularization recovered by the discrepancy principle is $\mu=7.9\cdot 10^{-3}$. The relative reconstruction error associated to isotropic Tikhonov regularization is $0.1902$. The local orientation parameter recovered by the new bilevel optimization method are displayed in Figure \ref{fig: tomo_theta}. The history of relevant quantities that monitor the progress of the new bilevel optimization approach to bilevel optimization are displayed in Figure \ref{fig: tomo_vals}. 

\subsection{Sparse Dix velocity inversion}
Sparse Dix velocity inversion is a seismic application that can significantly benefit from the proposed anisotropic regularization. This process involves estimating subsurface interval seismic velocity from root-mean-square (RMS) velocities. For a horizontally layered earth model, the forward operator is a causal integration matrix (a lower triangular matrix of ones) that links each column of the interval velocity matrix to the corresponding column in the RMS velocity matrix. In practice, due to computational limitations, the velocity analysis of common-depth-point (CDP) gathers provides estimates of RMS velocities at sparse locations. However, we need both the RMS velocity and the interval velocity at dense spatial grids. Sparse Dix velocity inversion addresses this by simultaneously performing the inversion from RMS velocity to interval velocity and interpolating to all CDP locations. This procedure is severely ill-conditioned, necessitating proper regularization to stabilize the solution \cite{Gholami_2019_3DD}.

\begin{center}
\begin{figure}
\begin{tabular}{cccc}
\hspace{-1.3cm}\small{(a)} & 
\hspace{-1.3cm}\small{(b)} & 
\hspace{-1.3cm}\small{(c)} & 
\hspace{-1.3cm}\small{(d)} \vspace{-0.1cm}\\
\hspace{-1.3cm}\includegraphics[width = 5cm]{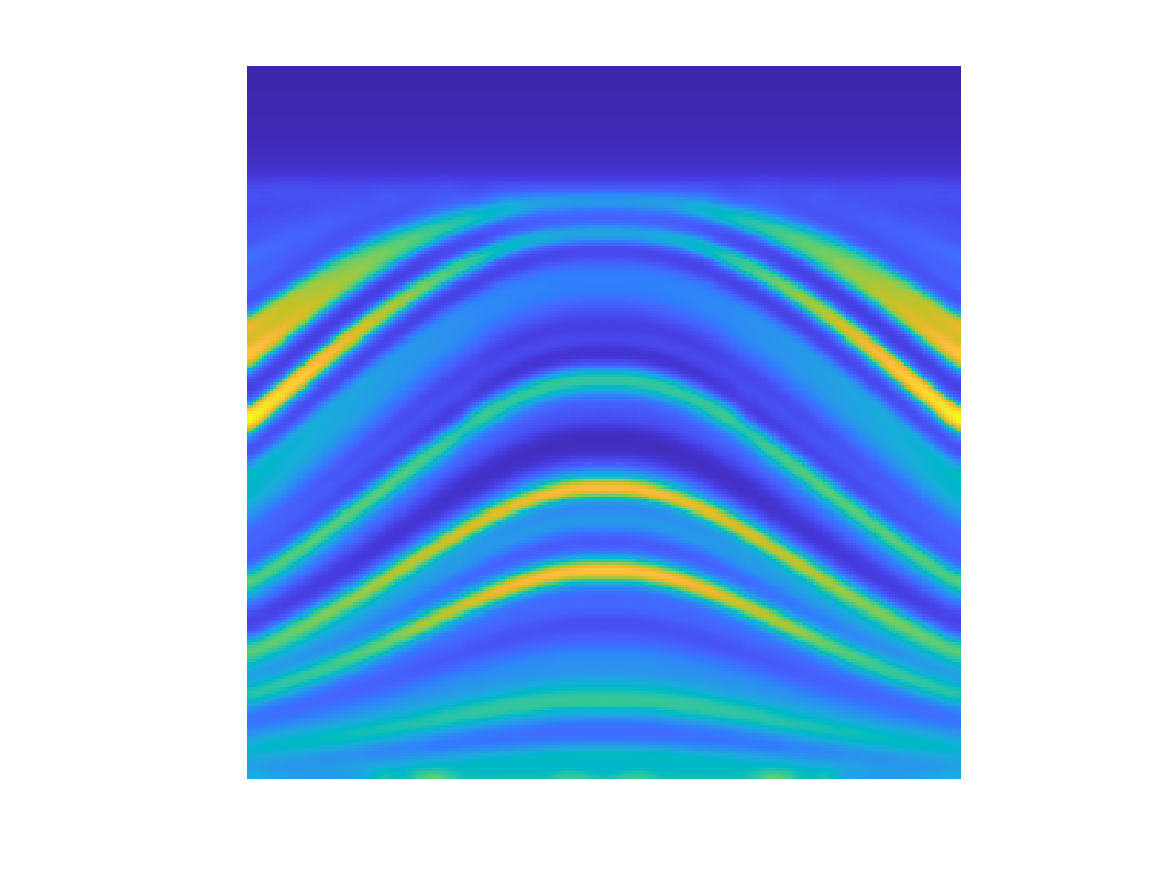} &
\hspace{-1.3cm}\includegraphics[width = 5cm]{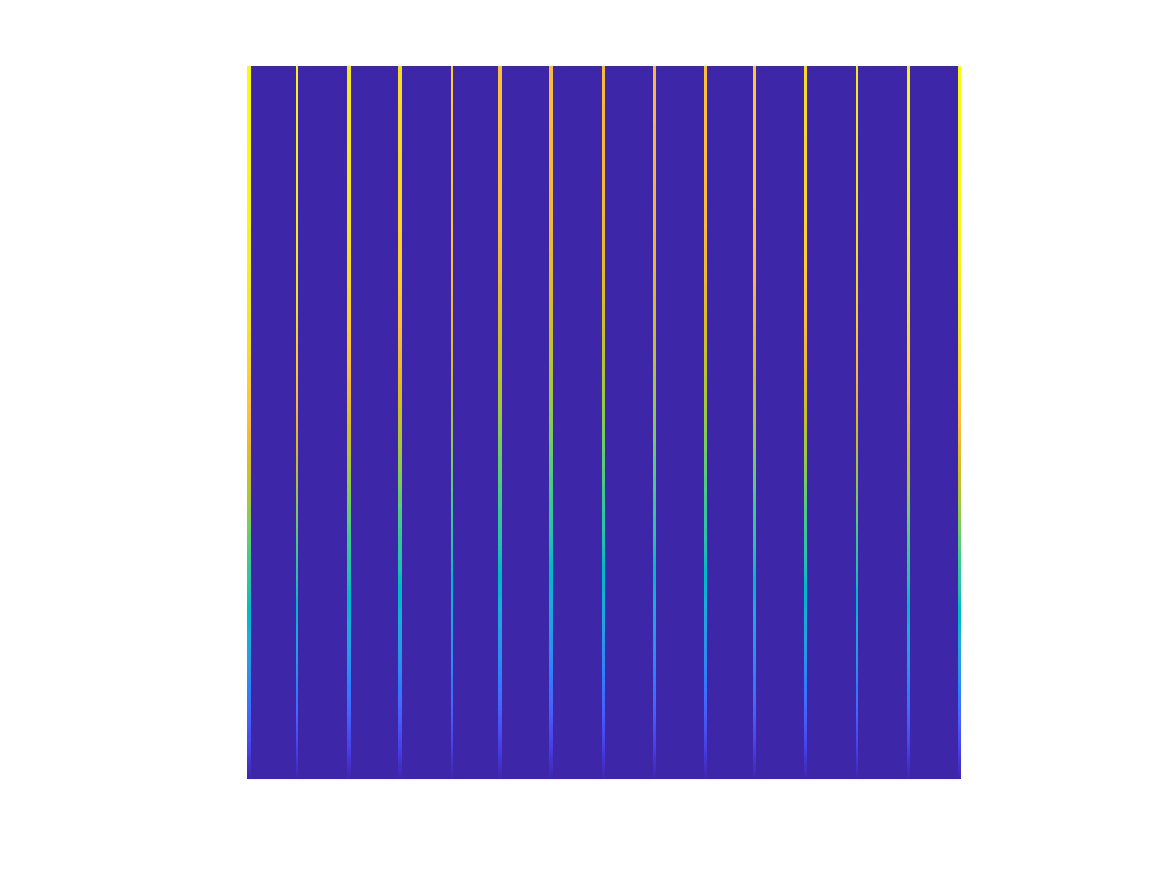} & 
\hspace{-1.3cm}\includegraphics[width = 5cm]{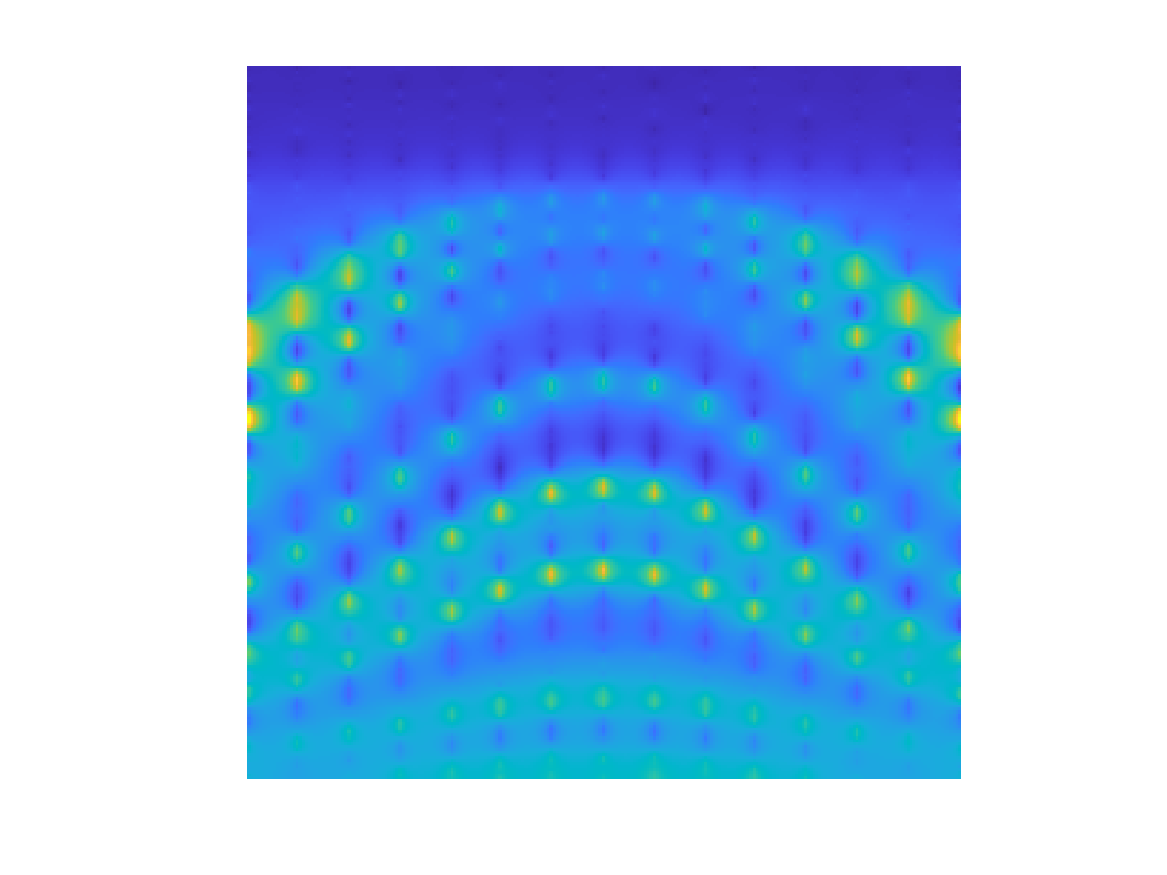} &
\hspace{-1.3cm}\includegraphics[width = 5cm]{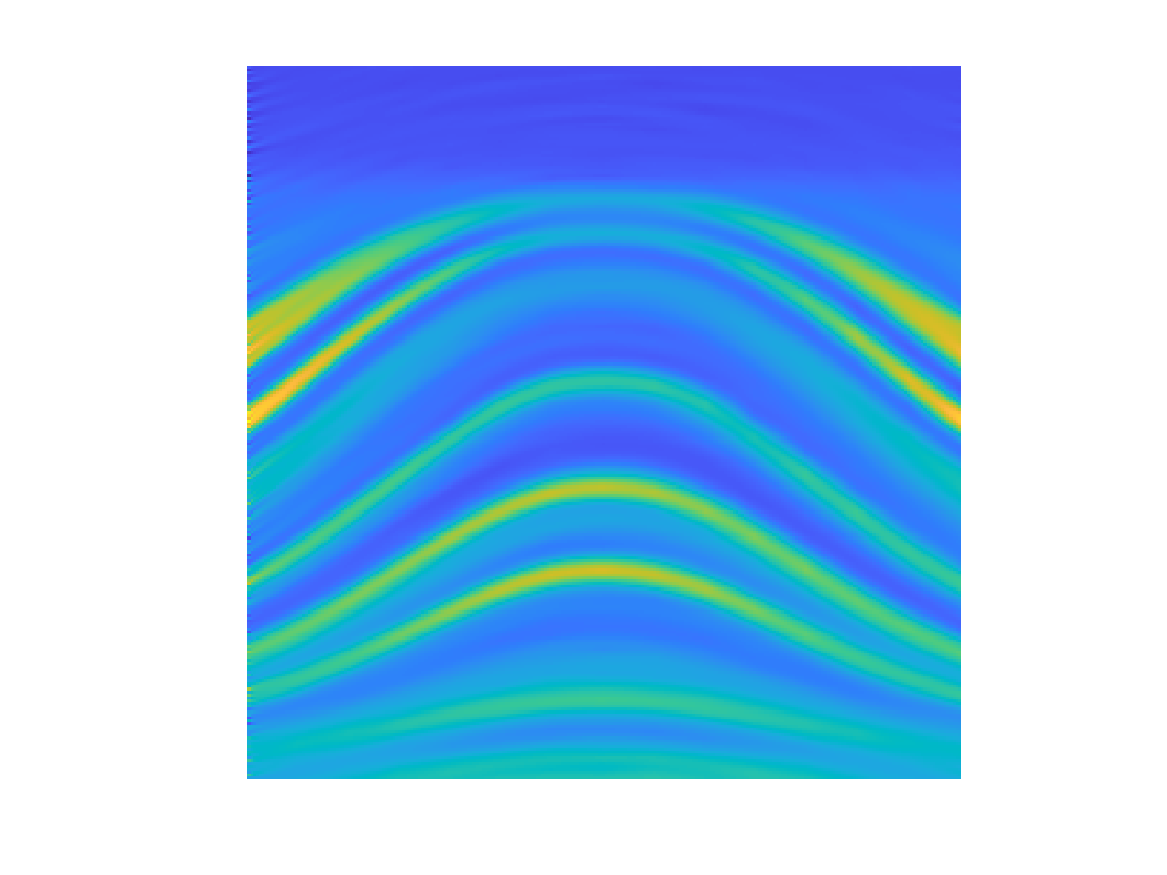}
\end{tabular}
\caption{Dix interpolation test problem. (a) exact parameters; (b) available data; (c) parameters recovered by the (isotropic) Tikhonov regularization method; (d) parameters recovered by the new anisotropic Tikhonov regularization method.}
\label{fig: dix_imgs}
\end{figure}
\end{center}

\begin{center}
\begin{figure}
\begin{tabular}{cc}
\hspace{-1.3cm}\small{(a)} & 
\hspace{-0.3cm}\small{(b)}\\
\hspace{-1.3cm}\includegraphics[height = 5cm]{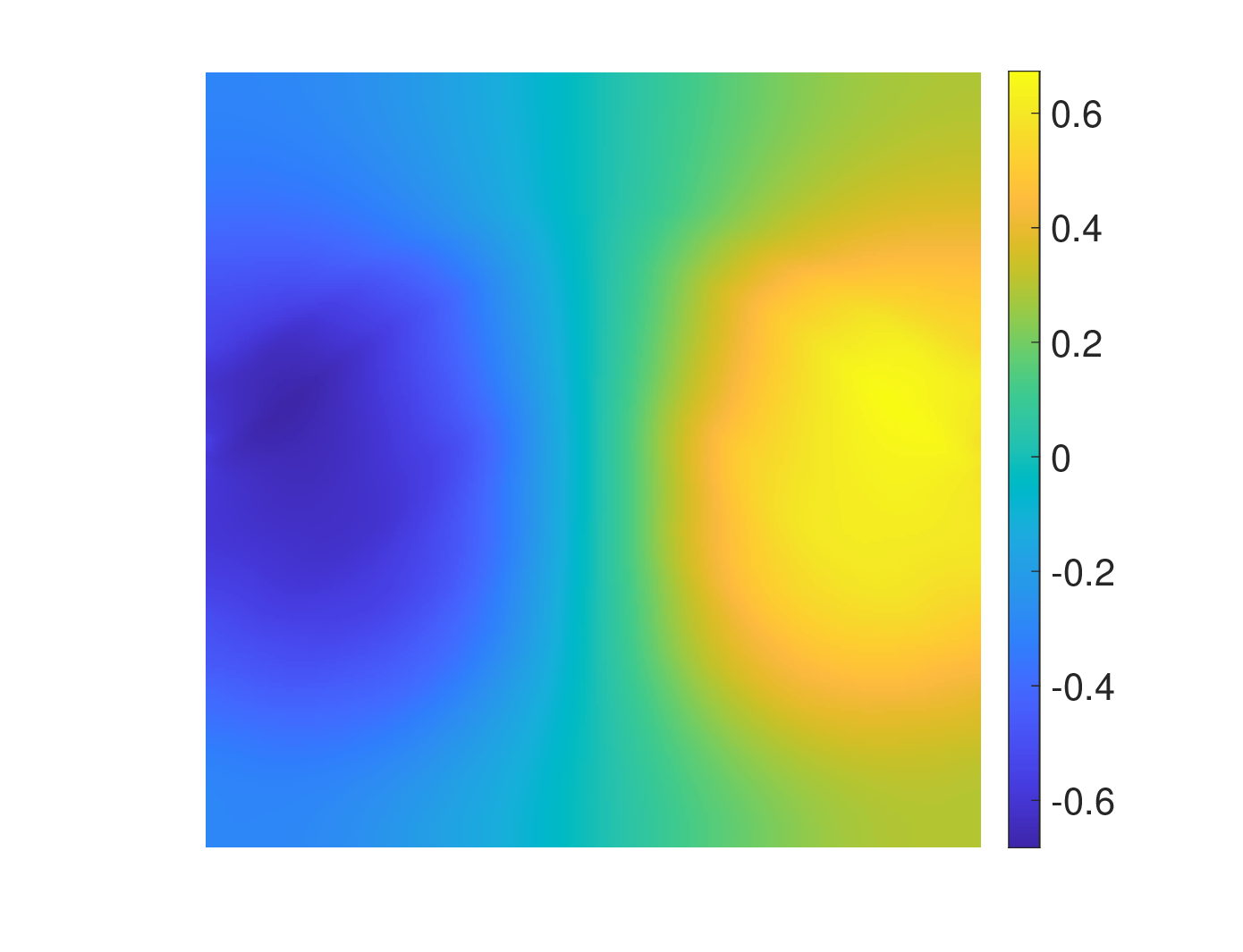} &
\hspace{-0.3cm}\includegraphics[height = 5cm]{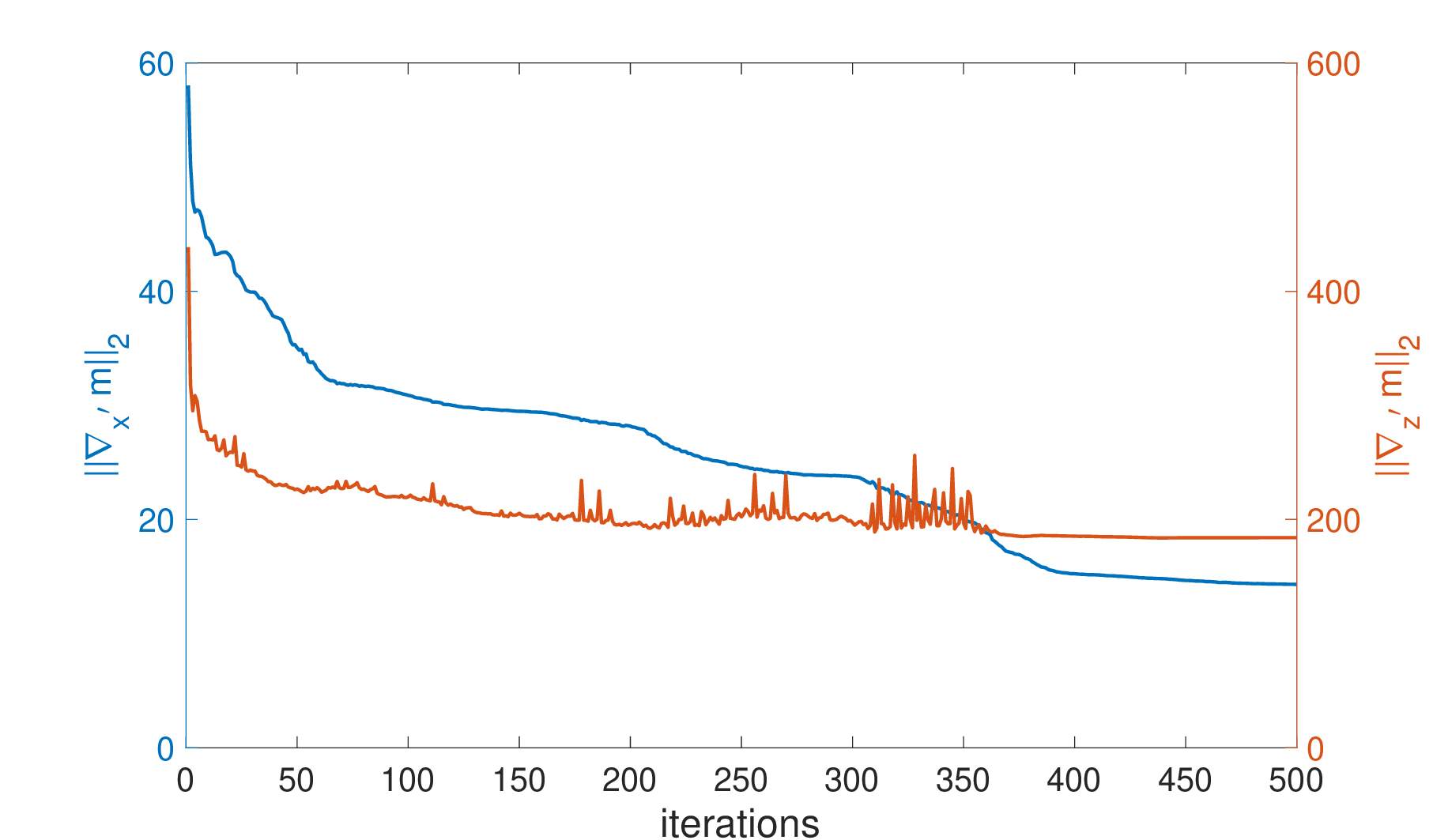}
\end{tabular}
\caption{Dix interpolation test problem. (a) pixel-wise orientation parameters recovered solving the bilevel optimization method \eqref{eq:bilevel}; (b) 2-norm of the directional derivatives along $x'$ and $z'$ versus L-BFGS-B iterations.}
\label{fig: dix_theta}
\end{figure}
\end{center}

\begin{center}
\begin{figure}
\begin{tabular}{ccc}
\hspace{-0.5cm}\small{(a)} & 
\hspace{-0.5cm}\small{(b)} &
\hspace{-0.5cm}\small{(c)}\\
\hspace{-0.5cm}\includegraphics[width = 5cm]{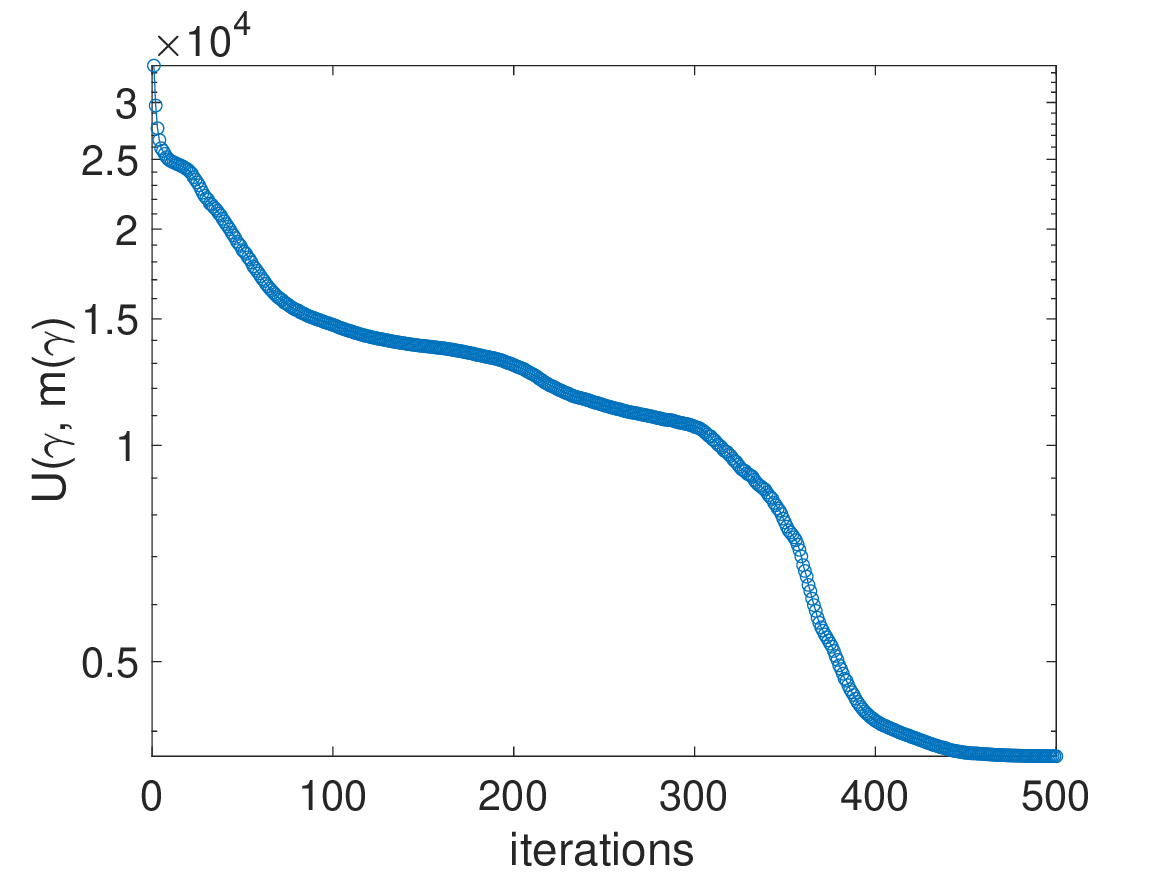} &
\hspace{-0.5cm}\includegraphics[width = 5cm]{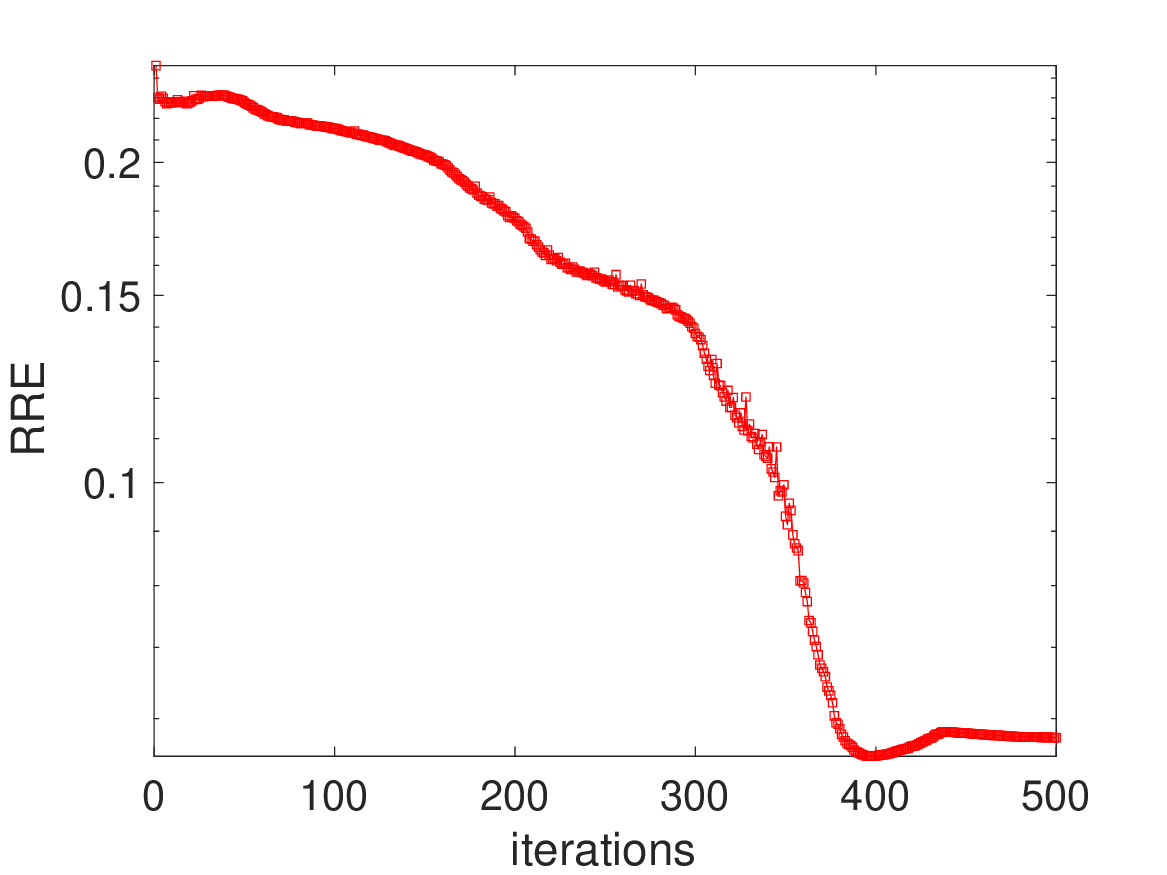} &
\hspace{-0.5cm}\includegraphics[width = 5cm]{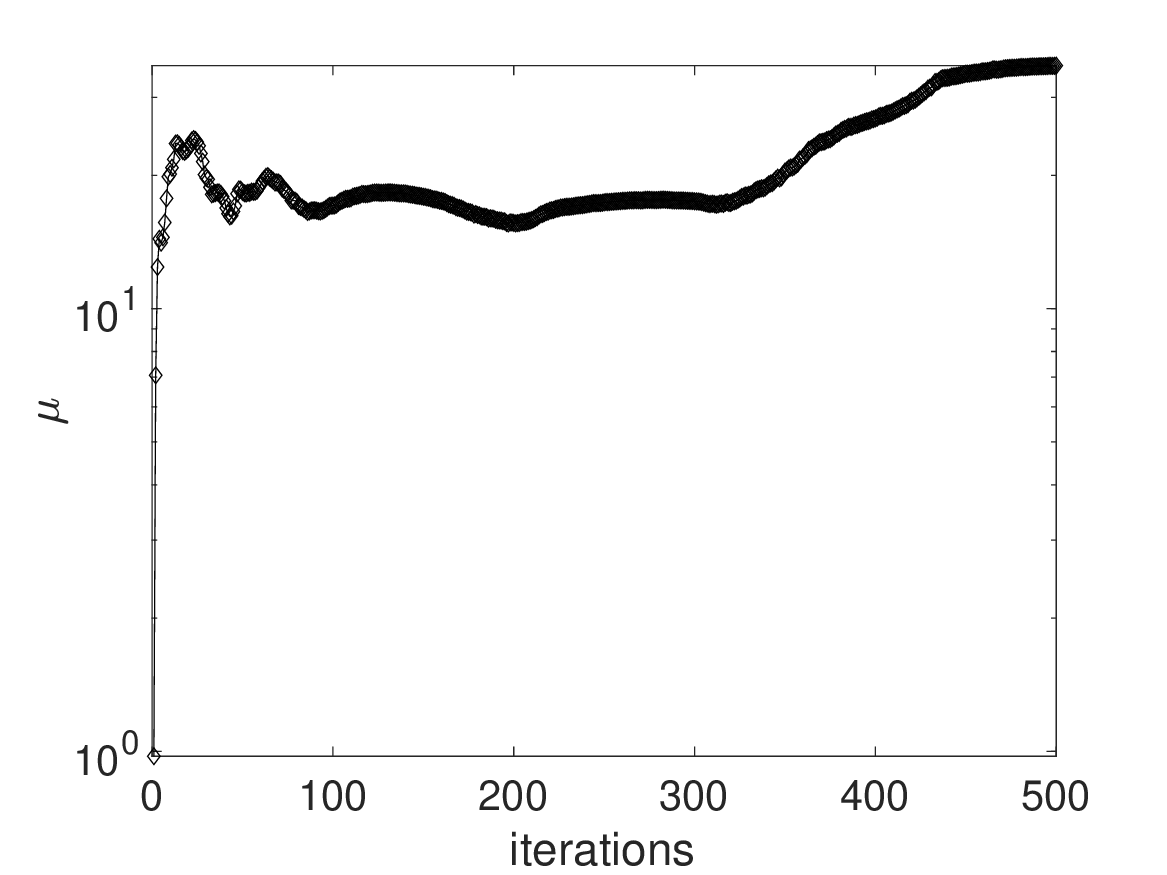}
\end{tabular}
\caption{Dix interpolation test problem. (a) upper level objective function values versus L-BFGS-B iterations; (b) relative reconstruction errors versus L-BFGS-B iterations; (c) Tikhonov regularization parameter $\mu$ versus L-BFGS-B iterations.}
\label{fig: dix_vals}
\end{figure}
\end{center}

The model parameters we aim to estimate are represented as a $250\times 250$ pixel image (Figure \ref{fig: dix_imgs}a). We assume that only $6\%$ of the CDP gathers are processed. The resulting RMS velocities at the selected CDP locations are shown in Figure \ref{fig: dix_imgs}b, which
 are affected by some Gaussian noise $\be$ of level $\|\be\|_2/\|\bG\bm_{\rm true}\|_2$ equal to $1.2\cdot 10^{-3}$. As weights for the gradient components we take $\sigma^{x'}=1$ and $\sigma^{z'} = 10^{-3}$. The parameters appearing in the upper level objective function \eqref{ULobjfcn} are taken as $\alpha=1$, $\beta=2\cdot 10^{3}$. The reconstructions using (isotropic) Tikhonov regularization and the new local anisotropic regularization strategy are shown in Figure \ref{fig: dix_imgs}. The regularization parameter for (isotropic) Tikhonov regularization recovered by the discrepancy principle is $\mu=1.6039$. The relative reconstruction error associated to the former is $0.2401$. The local orientation parameter recovered by the new bilevel optimization method are displayed in Figure \ref{fig: dix_theta}. The history of relevant quantities that monitor the progress of the new bilevel optimization approach to bilevel optimization are displayed in Figure \ref{fig: dix_vals}. 

\section{Conclusions}\label{sect:end}

A number of extensions of the proposed approach are possible. For instance, one may devise more efficient ways of handling large-scale computations, possibly by introducing preconditioners when inverting the Hessian of the lower level objective function. It would be interesting to allow space-variant contributions of the regularization by, e.g., adapting the weights in the weighted 2-norm to each spatial location. 
Also extension to local anisotropic regularizers expressed in the 1-norm would be meaningful. The main open questions that still remain are: (1) the automatic tuning of the regularization parameters appearing in the upper level functional and, (2) the introduction of a stopping criterion for the upper level iterations.

\bibliographystyle{plain}
\bibliography{biblio.bib}

\end{document}